\documentclass{article}
\usepackage{amssymb}
\usepackage{diagbox}
\usepackage[english]{babel}
\usepackage[pdftex]{graphicx}
\usepackage[utf8]{inputenc}
\usepackage{amsmath} 
\usepackage[T1]{fontenc}

\title{- Generalization of divergences by application of the deformed logarithm.
 -Applications to linear inverse problems
- Inversion algorithms.
- Version 1. -}
\author{ Henri Lantéri}

\begin{document}
\maketitle
\tableofcontents

\section{Reminders on the context.}
This work is related to the general context of inverse problems \cite{bertero1998}, \cite{Idier01a}.\\
The resolution of this type of problem implies to minimize (possibly under constraints), a function of discrepancy (Divergence) $D(p|q)$ between the measurements "$p$" and a physical model "$q$" of the considered phenomenon; the unknowns of the problem being the parameters of the model.\\
In the context of image deconvolution problems, for example, the constraints considered are typically the non-negativity constraint of the values of the intensities of the pixels of the reconstructed image, and the constraint of the sum of the total intensity of the reconstructed image. 
This approach implies on the one hand the definition of a deviation function and on the other hand the implementation of an algorithmic method of minimization under constraints of such a function.
These two aspects have been developed in a previous work\cite{lanteri2019}, \cite{lanteri2020}.\\
In order to handle in a simple way the constraint of sum of the unknown parameters, we have introduced in these previous works, the divergences invariant by change of scale on the unknowns (scale invariant divergences), and the notion of invariance factor.\\ 
 
\section{Objective of the present analysis.}
The purpose of this note is to extend the divergences analyzed in a previous work \cite{lanteri2019}, \cite{lanteri2020}, by application of the Deformed Logarithm in its most general form.\\ 
In a study on entropic divergences \cite{lanteri2021} \cite{lanteri2022}, we have analyzed the different forms of the deformed logarithm and their applications to this type of divergence.\\
A very general form of the deformed logarithm has been highlighted; it allows us to group the different expressions of the deformed logarithm under a unique form which includes all the others.\\
Beyond the applications linked to the divergences based on entropy, we propose here to extend the classical divergences by application of the generalized deformed Logarithm.

\section{A brief recall on the deformed logarithm.}
The various expressions of the deformed logarithm found in the literature have been presented in a previous work \cite{lanteri2021} \cite{lanteri2022} and the main properties of this function are given in \cite{naudts2002}.\\
We recall here the approach allowing to obtain the most general form of the deformed logarithm that we will use in the following.\\
In the context of statistical physics, by analogy with Shannon's entropy \cite{shannon1948}, a general expression of the entropy is written in the form:
\begin{equation}
S\left(p\right)=-\sum_{i}p_{i}\Lambda\left(p_{i}\right)  
\end{equation}
In this expression, the function $\Lambda\left(p_{i}\right)$ is the deformed logarithm or the natural logarithm in the case of the Shannon entropy which is expressed as:
\begin{equation}
S\left(p\right)=-\sum_{i}p_{i}\log p_{i}  
\end{equation}
This last expression can be obtained by using the definition proposed by Abe \cite{abe1997}:
\begin{equation}
	S(p)=\left[-\frac{d}{d\alpha}\sum_{i}p^{\alpha}_{i}\right]_{\alpha=1}\equiv \left[-\frac{d}{d\alpha}f(\alpha)\right]_{\alpha=1}
	\label{defabe}
\end{equation}
In this expression the involved derivative is the classical derivative.\\
From (\ref{defabe}), it is proposed to deform the classical derivative and to replace it by the Jackson ``$q$'' differential  \cite{jackson1909} \cite{jackson1910} which is expressed as:
 \begin{equation}
	\frac{df(\alpha)}{d(\alpha:t)}=\frac{f(\alpha t)-f(\alpha)}{\alpha t-\alpha} 
 \end{equation}

With this definition, we obtain:
\begin{equation}
	S(p)=\left[-\frac{\sum_{i}p^{\alpha t}_{i}-\sum_{i}p^{\alpha}_{i}}{\alpha t-\alpha}\right]_{\alpha=1}
\end{equation}
This corresponds to the Tsallis entropy \cite{tsallis1988}:
\begin{equation}
	S_{T}(p)=-\frac{\sum_{i}p^{t}_{i}-\sum_{i}p_{i}}{t-1}=-\frac{\sum_{i}p^{t}_{i}-1}{t-1}=-\sum_{i}p_{i}\left(\frac{p^{t-1}_{i}-1}{t-1}\right) 
\end{equation}
The corresponding deformed logarithm is then:
\begin{equation}
	Log_{T}(x)=\frac{x^{t-1}-1}{t-1}
\label{logdtsa}
\end{equation}
An alternative is, still from (\ref{defabe}), to adopt another definition of the differential given by Mc Anally \cite{mcanally1995} which is expressed:
\begin{equation}
	\frac{df(\alpha)}{d(\alpha:t)}=\frac{f(\alpha t)-f(\alpha t^{-1})}{\alpha t-\alpha t^{-1}}
\end{equation}
We thus have an invariance $t\leftrightarrow t^{-1}$ and we obtain:
\begin{equation}
	S(p)=\left[-\frac{\sum_{i}p^{\alpha t}_{i}-\sum_{i}p^{\alpha t^{-1}}_{i}}{\alpha t-\alpha t^{-1}}\right]_{\alpha=1}=-\frac{\sum_{i}p^{t}_{i}-\sum_{i}p^{t^{-1}}_{i}}{t- t^{-1}}
\end{equation}
We thus recover the entropy of Abe \cite{abe1997}; the corresponding deformed logarithm is expressed:
\begin{equation}
	Log_{A}(x)=\frac{x^{t-1}-x^{t^{-1}-1}}{t-t^{-1}}
\end{equation}
Another form of generalized 2-parameter differential, originally used by Chakrabarti and Jagannathan \cite{chakrabarti1991}, is proposed by Borges and Roditi \cite{borges1998}; it consists in using the following definition:
\begin{equation}
	\frac{d_{ab}f(\alpha)}{d_{ab}(\alpha)}=\frac{f(a \alpha)-f(b \alpha)}{a \alpha-b \alpha}
\end{equation}
This formulation extends the definitions previously proposed, it leads, with the definition (\ref{defabe}) to the expression of the entropy:
\begin{equation}
\begin{split}
	S_{ab}(p)=\left[-\frac{\sum_{i}p^{a \alpha}_{i}-\sum_{i}p^{b \alpha}_{i}}{a \alpha-b \alpha}\right]_{\alpha=1}&=-\frac{\sum_{i}p^{a}_{i}-\sum_{i}p^{b}_{i}}{a-b}\\&=-\sum_{i}p_{i}\left( \frac{p^{a-1}_{i}-p^{b-1}_{i}}{a-b}\right) 
	\label{entgen}
\end{split}
\end{equation}
This very general form of entropy (and the associated deformed logarithm) is mentioned by Wada and Scarfone \cite{wada2010}; it is based on the earlier work of Sharma and Taneja \cite{sharma1975} and Mittal \cite{mittal1975} and allows of course to find the expressions cited above.\
The range of values of the parameters ``$a$'' and ``$b$'' given in \cite{furuichi2010} and \cite{borges1998} are the following:
\begin{equation}
	0< a\leq 1\leq b
	\label{cond1}
\end{equation}
or:
\begin{equation}
 0<b\leq 1\leq a
\label{cond2}
\end{equation}
The expression of the deformed logarithm deduced from this form of entropy is given by:
\begin{equation}
	Log_{ab}(x)\equiv Log_{d}(x)=\frac{x^{a-1}-x^{b-1}}{a-b}
	\label{logdéformé}
\end{equation}
The deformed logarithm thus defined has the concavity properties of the classical logarithm in a domain of values of the parameters mentionned above.\\
Note that taking into account (\ref{cond1}) and (\ref{cond2}), we have:
\begin{equation}
\frac{a-1}{a-b}>0\;\;;\;\;\frac{b-1}{a-b}<0
\end{equation}
The special cases developed in \cite{lanteri2021} can be found by making the following adaptations:\\
\begin{itemize}
\item  Shannon entropy \cite{shannon1948}: $a\rightarrow 1$ et  $b\rightarrow 1$ taking into account (\ref{cond1}) ou (\ref{cond2})\\
\item Tsallis entropy \cite{tsallis1988} $a=t,$\ \ \ $  b=1$\\
\item Kaniadakis entropy \cite{kaniadakis2002}: $a=1+K,$\ \ \ $ b=1-K$\\
\item Abe entropy\cite{abe1997}: $a=z,$\ \ \ $ b=\frac{1}{z}$\\
\item ``$\gamma$'' entropy \cite{kaniadakis2005}: $a=2\gamma+1,$\ \ \ $ b=1-\gamma$\\
\item 2 parameters entropy (KLS) \cite{kaniadakis2005}: $a=1+r+K,$\ \ \ $ b=1+r-K$\\
\end{itemize}

In addition, while for the Natural Logarithm we have:
\begin{equation}
	Log(x.y)=Log(x)+Log(y)
	\label{lognprod}
\end{equation}
, for the deformed logarithm, this relation is not true anymore, and we have:
\begin{equation}
	Log_{d}(x.y)\neq Log_{d}(x)+Log_{d}(y)
	\label{logdprod}
\end{equation}
Note that the natural logarithm is obtained from the relation (\ref{logdéformé}) by taking for example $a=1+\epsilon$, $b=1-\epsilon$ and by performing the passage to the limit $\epsilon\rightarrow 0$.\\
The details of the derivation allowing this transition are given in Appendix 1.\\
The relations (\ref{lognprod}) and (\ref{logdprod}) imply different processes depending on whether the logarithmic divergences are obtained by applying the natural logarithm or the deformed logarithm.\\
We highlight in the subsequent developments the differences between the two approaches.

\section{Recalls on the scale invariant divergences - Invariance factor.}
This notion has been extensively presented in \cite{lanteri2019} and \cite{lanteri2020}, we only remind here the main points.\\
Starting from a divergence $D(p\|q)$, we propose to build a divergence $D(p\|Kq)=DI(p\|q)$ ( where "$K(p,q)$" is a positive scalar), such that $DI(p\|q)$ remains invariant when "$q$" is multiplied by a positive constant.\\
1- For a divergence $D(p|q)$, the nominal invariance factor $K_{0}(p,q)$ specific to the considered divergence is obtained as a solution of the equation:
\begin{equation}
	\sum_{i}\frac{\partial D(p_{i}\|Kq_{i})}{\partial K}=0
\end{equation}
if an explicit solution is available.\\

2- The general properties of invariance factors are as follows:\\
* the invariance factor $K(p,q)$ is a positive scalar.\\
* The components of the vector $\left[K(p,q).q\right]$ are invariant if "$q$" is multiplied by a constant positive factor.\\
* All factors $K(p,q)$ having the above properties make any divergence scale invariant.\\
They are solutions of the differential equation:	
\begin{equation}
	K(p,q)+\sum_{j}q_{j}\frac{\partial K(p,q)}{\partial q_{j}}=0	
\end{equation}
The nominal invariance factors are of course solutions of the above differential equation.\\

3- Whether the invariance factor is the nominal invariance factor or not, we have the following fundamental property:
\begin{equation}
\sum_{j}q_{j}\frac{\partial D(p\|K(p,q).q)}{\partial q_{j}}=0	
\end{equation}

\section{General forms of the divergences and application of the deformed logarithm.}
\subsection{Classical forms of the divergences.}
A divergence is by definition a positive quantity.
As we have indicated in section l, it expresses a difference between 2 data fields: the "$p$" field which for us will always represent the measurements of the considered phenomenon, and the "$q$" field which will represent the physical model of the same phenomenon.\\
The parameters of the model are the unknowns of the problem.\\
This being so, many divergences, whether they are non-invariant or whether they are made invariant by the introduction of an invariance factor whatever it is, can be expressed in the form of a difference of 2 positive terms.\\

\textbf{*Remark}: When the nominal invariance factor has no explicit expression, one will be led to use systematically, as invariance factor, a specific expression; the scale invariant divergences thus obtained will be discussed on a case by case basis.\\

These particular situations being put aside (for the moment), a very large number of divergences can be written in the following general form::
\begin{equation}
 D(p\|q)=\left[A\left( p,q\right)-B\left( p,q\right)\right]\;\;;\;\;A\left( p,q\right)>0\;\;;\;\;B\left( p,q\right)>0  
\label{divD}	
\end{equation}

From this most common relation, 2 cases appear: either $B=X+Y$, or $B=X.Y$ which leads to the classical expressions:
\begin{equation}
D1(p\|q)=\left\lbrace A\left(p,q\right)-\left[X\left(p,q\right)+Y\left(p,q\right)\right] \right\rbrace  
\label{divD1}	
\end{equation}

\begin{equation}
D2(p\|q)=\left\lbrace A\left(p,q\right)-\left[X\left(p,q\right).\,Y\left(p,q\right)\right] \right\rbrace  
\label{divD2}	
\end{equation}

Moreover, if the expressions $A,\:B,\:X,\:Y$ are separable, we have:

\begin{equation}
\begin{split}
&A\left(p,q\right)=\sum_{i}A\left(p_{i},q_{i}\right)=\sum_{i}A_{i}\ ;\ \  A_{i}>0\ \ \forall i\\&
B\left(p,q\right)=\sum_{i}B\left(p_{i},q_{i}\right)=\sum_{i}B_{i}\ ;\ \  B_{i}>0\ \ \forall i\\&
X\left(p,q\right)=\sum_{i}X\left(p_{i},q_{i}\right)=\sum_{i}X_{i}\ ;\ \  X_{i}>0\ \ \forall i\\&
Y\left(p,q\right)=\sum_{i}Y\left(p_{i},q_{i}\right)=\sum_{i}Y_{i}\ ;\ \  Y_{i}>0\ \ \forall i
\end{split}	
\end{equation}

\textbf{When applying deformed logarithm on such expressions, insofar as this transformation must preserve the positivity of the divergence obtained, the precautions to be taken are the following:\\
1 - A term of the form $\left[X\left( p,q\right)+Y\left( p,q\right)\right]$ must be treated as an indissociable block (which seems obvious).\\
2 - Given (\ref{logdprod}), a term of the form $\left[X\left( p,q\right).\:Y\left( p,q\right)\right]$ must be regarded as an indissociable block.}\\

For $D$ type divergences, the logarithmic forms proposed in the literature are obtained by applying the natural or deformed logarithm on each term of the divergence which leads with the natural Logarithm, to:
\begin{equation}
LD(p\|q)=\left[\log A\left(p,q\right)-\log B\left(p,q\right)\right]	
\end{equation}
or, with the deformed Logarithm, to:
\begin{equation}
L_{d}D(p\|q)=\left[\log_{d} A\left(p,q\right)-\log_{d} B\left(p,q\right)\right]
\label{logdD}		
\end{equation}
For the form D1, this operation does not cause any particular problem;
by applying the natural logarithm on the two terms of the divergence, $A$ and $X+Y$, it comes:
\begin{equation}
LD1(p\|q)=\left\lbrace \log A\left(p,q\right)-\log \left[ X\left(p,q\right)+ Y\left(p,q\right)\right] \right\rbrace 
\label{lognatD1}	
\end{equation}
Similarly, applying the deformed logarithm, we will have:
\begin{equation}
L_{d}D1(p\|q)=\left\lbrace \log_{d} A\left(p,q\right)-\log_{d} \left[ X\left(p,q\right)+ Y\left(p,q\right)\right] \right\rbrace 
\label{logdD1}	
\end{equation}

The case of the divergence $D2(p\|q)$ leads by application of the Natural Logarithm to:
\begin{equation}
LD2(p\|q)=\left[\log A\left(p,q\right)-\log X\left(p,q\right)-\log Y\left(p,q\right)\right]
\label{lognatD2}		
\end{equation}
On the contrary, with the deformed logarithm, the decomposition of the product term is impossible given (\ref{logdprod}), and we will have:
\begin{equation}
L_{d}D2(p\|q)=\left\lbrace\log_{d} A\left(p,q\right)-\log_{d}\left[ X\left(p,q\right). Y\left(p,q\right)\right]\right\rbrace 
\label{logdD2}		
\end{equation}

\subsection{Particular forms of divergences.}

\subsubsection{Divergences involving the natural logarithm.}
Some divergences proposed in the literature are based on the natural Logarithm function, it is for example the well known case of the Kullback-Leibler \cite{kullback1951} divergence based on the Shannon entropy and of all the divergences deduced from it.\\
In this case, the deformed logarithm is simply introduced in replacement of the natural logarithm.\\
This extension has been considered in a previous work \cite{lanteri2021} dealing specifically with entropic divergences and can be implemented systematically..\\

\textbf{* WARNING}: Some divergences, for instance the "Gamma divergence" \cite{cichocki2010} are obtained by applying the natural Logarithm function on divergences of the form $D2(p|q)$ (\ref{divD2}) according to the expression (\ref{lognatD2}).\\
For such divergences, it is of course incorrect to replace the natural logarithm by the deformed logarithm because this would be equivalent to considering that the property (\ref{lognprod}) specific to the natural logarithm is also valid for the deformed logarithm (which is incorrect).

\subsubsection{Specific forms of scale invariant divergences.}
As previously mentioned, for some divergences, it is impossible to derive an explicit expression for the nominal invariance factor.\\
In these particular cases, one will systematically use the invariance factor $K^{*}=\frac{\sum_{j}p_{j}}{\sum_{j}q_{j}}$  which is the nominal invariance factor for the Kullback-Leibler divergence.\\
When such an invariance factor is used, except for a multiplicative factor, which depends only on $\sum_{j}p_{j}$, the invariant divergence obtained has the same expression as the initial divergence provided that we replace the variables "$p$" and "$q$" by normalized variables $\overline{p}_{i}=\frac{p_{i}}{\sum_{j}p_{j}}$ and $\overline{q}_{i}=\frac{q_{i}}{\sum_{j}q_{j}}$ . \\
This being so, everything in Section (5.1) continues to hold with the new variables.

\subsection{Expressions of Divergences and of their Gradients with respect to "$q$".}

\subsubsection{Divergence of type LD.}
By applying the generalized logarithm given by (\ref{logdéformé}), the divergence (\ref{divD}) becomes:
\begin{equation}
L_{d}D(p\|q)=\left\{\frac{A^{a-1}\left(p,q\right)-A^{b-1}\left(p,q\right)}{a-b}-\frac{B^{a-1}\left(p,q\right)-B^{b-1}\left(p,q\right)}{a-b}\right\}
\label{logdDexplicite}
\end{equation}
Its gradient with respect to "$q$" is expressed as, $\forall j$:
\begin{equation}
\begin{split}
	\frac{\partial L_{d}D(p\|q)}{\partial q_{j}}=&\left\{\left[\frac{a-1}{a-b}A^{a-2}-\frac{b-1}{a-b}A^{b-2}\right]\frac{\partial A}{\partial q_{j}}\right\}\\&-\left\{\left[\frac{a-1}{a-b}B^{a-2}-\frac{b-1}{a-b}B^{b-2}\right]\frac{\partial B}{\partial q_{j}}\right\}
	\label{gradlogdD}
\end{split}	
\end{equation}
This expression is further simplified when the basic divergences considered are separable, we have indeed:
\begin{equation}
\begin{split}	
&\frac{\partial A}{\partial q_{j}}=\frac{\partial A(p_{j},q_{j})}{\partial q_{j}}\equiv \frac{\partial Aj}{\partial q_{j}}
\\&\frac{\partial B}{\partial q_{j}}=\frac{\partial B(p_{j},q_{j})}{\partial q_{j}}\equiv \frac{\partial Bj}{\partial q_{j}}
\label{dérivsimpl}
\end{split}
\end{equation}

\subsubsection{Divergence of type LD1.}
By applying the generalized logarithm given by (\ref{logdéformé}), the divergence (\ref{divD1}) becomes:
\begin{equation}
\begin{split}
L_{d}D1(p\|q)=&\left[\frac{A^{a-1}\left(p,q\right)-A^{b-1}\left(p,q\right)}{a-b}\right]\\&-\left[\frac{(X+Y)^{a-1}\left(p,q\right)-(X+Y)^{b-1}\left(p,q\right)}{a-b}\right] 
\label{logdD1explicite}
\end{split}
\end{equation}
Its gradient with respect to "$q$" is expressed as, $\forall j$:
\begin{equation}
\begin{split}
	\frac{\partial L_{d}D1(p\|q)}{\partial q_{j}}&=\left\{\left[\frac{a-1}{a-b}A^{a-2}-\frac{b-1}{a-b}A^{b-2}\right]\frac{\partial A}{\partial q_{j}}\right\}\\&-\left\{\left[\frac{a-1}{a-b}(X+Y)^{a-2}-\frac{b-1}{a-b}(X+Y)^{b-2}\right]\frac{\partial (X+Y)}{\partial q_{j}}\right\}
	\label{gradlogdD1}
\end{split}	
\end{equation}
This expression is further simplified when the basic divergences considered are separable, we have indeed:
\begin{equation}
\frac{\partial A}{\partial q_{j}}=\frac{\partial A(p_{j},q_{j})}{\partial q_{j}}\:\:\:\:;\:\:\:\:\frac{\partial (X+Y)}{\partial q_{j}}=\frac{\partial X(p_{j},q_{j})}{\partial q_{j}}+\frac{\partial Y(p_{j},q_{j})}{\partial q_{j}}
\end{equation}

\subsubsection{Divergence of type LD2.}
Here, the problem becomes a little more complicated due to the property (or non-property) (\ref{logdprod}) of the deformed logarithm, indeed, the divergence we are considering was written (\ref{logdD2}):
\begin{equation}
L_{d}D2(p\|q)=\left\{\log_{d} A\left(p,q\right)-\log_{d} \left[X\left(p,q\right). Y\left(p,q\right)\right]\right\}
\label{logdD2bis}	
\end{equation}
If we explicit the deformed logarithm function, we get:
\begin{equation}
L_{d}D2(p\|q)=\left\{\frac{A^{a-1}-A^{b-1}}{a-b}-\left[\frac{(X.Y)^{a-1}-(X.Y)^{b-1}}{a-b}\right]\right\}
\label{logdD2bisdev}		
\end{equation}
Its gradient with respect to "$q$" is expressed as, $\forall j$:
\begin{equation}
\begin{split}
	\frac{\partial L_{d}D2(p\|q)}{\partial q_{j}}=&\left\{\left[\frac{a-1}{a-b}A^{a-2}-\frac{b-1}{a-b}A^{b-2}\right]\frac{\partial A}{\partial q_{j}}\right\}\\&-\left\{\left[\frac{a-1}{a-b}(X.Y)^{a-2}-\frac{b-1}{a-b}(X.Y)^{b-2}\right]\frac{\partial (X.Y)}{\partial q_{j}}\right\}
	\label{gradlogdD2bis}	
\end{split}	
\end{equation}

In certain cases, taking into account the separability of the terms $A(p,q)$, $X(p,q)$ and $Y(p,q)$, we have:
\begin{equation}
\frac{\partial A(p,q)}{\partial q_{j}}=	\frac{\partial A(p_{j},q_{j})}{\partial q_{j}}=\frac{\partial A_{j}}{\partial q_{j}}	
\end{equation}
\begin{equation}
\frac{\partial (X.Y)}{\partial q_{j}}=X\frac{\partial Y(p_{j},q_{j})}{\partial q_{j}}+Y\frac{\partial X(p_{j},q_{j})}{\partial q_{j}}=X\frac{\partial Y_{j}}{\partial q_{j}}+Y\frac{\partial X_{j}}{\partial q_{j}}	
\end{equation}
In the case of the natural logarithm expressed by (\ref{lognatD2}) the gradient can be obtained without difficulty either by direct computation, or from the previous relations by following the procedure indicated in Appendix 1 and leads to:
\begin{equation}
	\frac{\partial LD2(p\|q)}{\partial q_{j}}=\left\lbrace\frac{1}{A}\frac{\partial A}{\partial q_{j}}-\left[\frac{1}{X}\frac{\partial X}{\partial q_{j}}+\frac{1}{Y}\frac{\partial Y}{\partial q_{j}}\right]\right\rbrace 
\end{equation}

\section{Application to $(\alpha)$, $(\beta)$ and $(\alpha \beta)$ divergences.}
\subsection{Alpha ($\alpha$) divergence.}
This divergence already analyzed in \cite{lanteri2019}, \cite{lanteri2020} and in the papers cited there \cite{amari2009},\cite{cichocki2010},\cite{cichocki2011} is written in the form:
\begin{equation}
	D_{\alpha}\left(p\|q\right)=\frac{1}{\alpha\left(\alpha-1\right)}\left\lbrace \sum_{i}p^{\alpha}_{i}q^{1-\alpha}_{i}-\left(\alpha\sum_{i} p_{i}+\left(1-\alpha\right)\sum_{i}q_{i}\right)\right\rbrace   
	\label{alphadiv}
\end{equation}
It can be expressed in the equivalent form:
\begin{equation}
	D_{\alpha}\left(p\|q\right)=\frac{1}{\alpha\left(\alpha-1\right)} \sum_{i}p^{\alpha}_{i}q^{1-\alpha}_{i}- \frac{1}{\alpha-1}\sum_{i} p_{i}+\frac{1}{\alpha}\sum_{i}q_{i}   
	\label{alphadivbis}
\end{equation}
The decomposition of this divergence under the form $A-B$ with $A>0$ and $B>0$ implies to distinguish 3 cases according to the values of the parameter "$\alpha$".\\\\

\textbf{*Case 1:}  $0<\alpha<1$, we then have:
\begin{equation}
A=\frac{1}{1-\alpha}\sum_{i} p_{i}+\frac{1}{\alpha}\sum_{i} q_{i}\;\;;\;\;B=\frac{1}{\alpha\left(1-\alpha\right)}\sum_{i}p^{\alpha}_{i}q^{1-\alpha}_{i}
\end{equation}
By applying the deformed logarithm on each term of the difference, we obtain ((\ref{logdDexplicite}) whose gradient with respect to "$q$" is given by ((\ref{gradlogdD}); taking into account the separability of the terms $A(p,q)$ and $B(p,q)$, we have:
\begin{equation}
\frac{\partial A}{\partial q_{j}}=\frac{\partial A_{j}}{\partial q_{j}}=\frac{1}{\alpha}\;\;;\;\;\frac{\partial B}{\partial q_{j}}=\frac{\partial B_{j}}{\partial q_{j}}=\frac{1}{\alpha}p^{\alpha}_{j}q^{-\alpha}_{j}
\end{equation}

The expression of the opposite of the gradient with respect to "$q$" is written, $\forall j$:
\begin{equation}
\begin{split}
-\frac{\partial L_{d}D_{\alpha}(p\|q)}{\partial q_{j}}=&-\frac{1}{\alpha}\left[\frac{a-1}{a-b}A^{a-2}-\frac{b-1}{a-b}A^{b-2}\right]\\&+
\frac{1}{\alpha}\left[\frac{a-1}{a-b}B^{a-2}-\frac{b-1}{a-b}B^{b-2}\right]p^{\alpha}_{j}q^{-\alpha}_{j}	
\end{split}	
\end{equation}
For an algorithmic use, the expression of the opposite of the gradient is written in the form $U-V$ with $U>0$ and $V>0$, consequently taking into account the conditions (\ref{cond1}) and (\ref{cond2}), it comes:
\begin{equation}
\begin{split}
&U_{j}=\frac{1}{\alpha}\left[\frac{a-1}{a-b}B^{a-2}-\frac{b-1}{a-b}B^{b-2}\right]p^{\alpha}_{j}q^{-\alpha}_{j}
\\&
V_{j}=\frac{1}{\alpha}\left[\frac{a-1}{a-b}A^{a-2}-\frac{b-1}{a-b}A^{b-2}\right]
\end{split}
\end{equation}

\textbf{*Case 2:}  $\alpha>1$, we have:
\begin{equation}
A=\frac{1}{\alpha}\sum_{i} q_{i}+\frac{1}{\alpha\left(\alpha-1\right)}\sum_{i}p^{\alpha}_{i}q^{1-\alpha}_{i}\;\;;\;\;B=\frac{1}{\alpha-1}\sum_{i} p_{i}
\end{equation}
By applying the deformed logarithm on each term of the difference, we obtain ((\ref{logdDexplicite}) whose gradient with respect to "$q$" is given by ((\ref{gradlogdD}); taking into account the separability of the terms $A(p,q)$ and $B(p,q)$, we have:
\begin{equation}
\frac{\partial A}{\partial q_{j}}=\frac{\partial A_{j}}{\partial q_{j}}=-\frac{1}{\alpha}p^{\alpha}_{j}q^{-\alpha}_{j}+\frac{1}{\alpha}\;\;;\;\;\frac{\partial B}{\partial q_{j}}=\frac{\partial B_{j}}{\partial q_{j}}=0
\end{equation}

The expression of the opposite of the gradient with respect to "$q$" is written, $\forall j$:
\begin{equation}
-\frac{\partial L_{d}D_{\alpha}(p\|q)}{\partial q_{j}}=\frac{1}{\alpha}\left[\frac{a-1}{a-b}A^{a-2}-\frac{b-1}{a-b}A^{b-2}\right]\left(p^{\alpha}_{j}q^{-\alpha}_{j}-1\right)
\end{equation}
For an algorithmic use, the expression of the opposite of the gradient is written in the form $U-V$ with $U>0$ and $V>0$, consequently taking into account the conditions (\ref{cond1}) and (\ref{cond2}), it comes:
\begin{equation}
\begin{split}
&U_{j}=\frac{1}{\alpha}\left[\frac{a-1}{a-b}A^{a-2}-\frac{b-1}{a-b}A^{b-2}\right]p^{\alpha}_{j}q^{-\alpha}_{j}
\\&
V_{j}=\frac{1}{\alpha}\left[\frac{a-1}{a-b}A^{a-2}-\frac{b-1}{a-b}A^{b-2}\right]
\end{split}
\end{equation}

\textbf{*Case 3:}  $\alpha<0$.\\
This case, which is a somewhat surprising one, must nevertheless be considered because it is part of the developments concerning the "$\alpha$" divergence.\\
We have then:
\begin{equation}
A=\frac{1}{1-\alpha}\sum_{i} p_{i}+\frac{1}{\alpha\left(\alpha-1\right)}\sum_{i}p^{\alpha}_{i}q^{1-\alpha}_{i}\;\;;\;\;B=-\frac{1}{\alpha}\sum_{i} q_{i}
\end{equation}
By applying the deformed logarithm on each term of the difference, we obtain ((\ref{logdDexplicite}) whose gradient with respect to "$q$" is given by ((\ref{gradlogdD}); taking into account the separability of the terms $A(p,q)$ and $B(p,q)$, we have:
\begin{equation}
\frac{\partial A}{\partial q_{j}}=\frac{\partial A_{j}}{\partial q_{j}}=-\frac{1}{\alpha}p^{\alpha}_{j}q^{-\alpha}_{j}\;\;;\;\;\frac{\partial B}{\partial q_{j}}=\frac{\partial B_{j}}{\partial q_{j}}=-\frac{1}{\alpha}
\end{equation}

The expression of the opposite of the gradient with respect to "$q$" is written, $\forall j$:
\begin{equation}
\begin{split}
-\frac{\partial L_{d}D_{\alpha}(p\|q)}{\partial q_{j}}=&\frac{1}{\alpha}\left[\frac{a-1}{a-b}A^{a-2}-\frac{b-1}{a-b}A^{b-2}\right]p^{\alpha}_{j}q^{-\alpha}_{j}\\&-\frac{1}{\alpha}\left[\frac{a-1}{a-b}B^{a-2}-\frac{b-1}{a-b}B^{b-2}\right]
\end{split}
\end{equation}
For an algorithmic use, the expression of the opposite of the gradient is written in the form $U-V$ with $U>0$ and $V>0$, consequently taking into account the conditions (\ref{cond1}) and (\ref{cond2}), it comes:
\begin{equation}
\begin{split}
&U_{j}=-\frac{1}{\alpha}\left[\frac{a-1}{a-b}B^{a-2}-\frac{b-1}{a-b}B^{b-2}\right]
\\&
V_{j}=-\frac{1}{\alpha}\left[\frac{a-1}{a-b}A^{a-2}-\frac{b-1}{a-b}A^{b-2}\right]p^{\alpha}_{j}q^{-\alpha}_{j}
\end{split}
\end{equation}

\subsection{Scale invariant Alpha ($\alpha$) divergence.}
This divergence has been studied in a previous work \cite{lanteri2019}, \cite{lanteri2020}.\\
By introducing its nominal invariance factor which is explicitly computable, and is expressed:
 \begin{equation}
 K_{0,\alpha}=\left[\frac{\sum_{i}p^{\alpha}_{i}q^{1-\alpha}_{i}}{\sum_{i}q_{i}}\right]^{\frac{1}{\alpha}}
 \end{equation}
 
We obtain the invariant form:
\begin{equation}
D_{\alpha}I\left(p\|q\right)=\underbrace{\frac{1}{1-\alpha}}_{T}\left[\underbrace{\sum_{i}p_{i}}_{A}-\underbrace{\left(\sum_{i}p^{\alpha}_{i}q^{1-\alpha}_{i}\right)^{\frac{1}{\alpha}}}_{X}\times\underbrace{\left(\sum_{i}q_{i}\right)^{1-\frac{1}{\alpha}}}_{Y}\right]
\label{alphadivinv}
\end{equation}
Disregarding the multiplicative factor "$T$", this divergence is of the form $D2$, moreover it is unnecessary to consider the different possible values of "$\alpha>0$" which all lead to identical results.\\\\
\textbf{Note: The case "$\alpha<0$" will be commented at the end of this paragraph.}\\\\
Since the terms "$X$" and "$Y$" are not separable, the calculation of the partial derivatives leads to:
\begin{equation}
\begin{split}
&\frac{\partial A}{\partial q_{j}}=0
\\&\frac{\partial X}{\partial q_{j}}=\frac{1-\alpha}{\alpha}\left(\sum_{i}p^{\alpha}_{i}q^{1-\alpha}_{i}\right)^{\frac{1}{\alpha}-1}p^{\alpha}_{j}q^{-\alpha}_{j}
\\&\frac{\partial Y}{\partial q_{j}}=\frac{\alpha-1}{\alpha}\left(\sum_{i}q_{i}\right)^{-\frac{1}{\alpha}}.\ 1
\label{derivpartiellesa}	
\end{split}	
\end{equation}
By applying the deformed logarithm, following the development of section (3.3.3), taking into account the term "$T$", we obtain the expression of the opposite of the gradient:
\begin{equation}
\begin{split}
	-\frac{\partial L_{d}D_{\alpha}I(p\|q)}{\partial q_{j}}=&-\frac{1}{\alpha}\left[\frac{a-1}{a-b}(X.Y)^{a-2}-\frac{b-1}{a-b}(X.Y)^{b-2}\right]\times\\&\left[\left(\frac{\sum_{i}p^{\alpha}_{i}q^{1-\alpha}_{i}}{\sum_{i}q_{i}}\right)^{\frac{1}{\alpha}}-\left(\frac{\sum_{i}p^{\alpha}_{i}q^{1-\alpha}_{i}}{\sum_{i}q_{i}}\right)^{\frac{1}{\alpha}-1}p^{\alpha}_{j}q^{-\alpha}_{j}\right]
\label{gradalphadivinv}
\end{split}	
\end{equation}
The decomposition of the opposite of the gradient in the form $U-V$ with $U>0$ and $V>0$ is expressed:
\begin{equation}
\begin{split}
&U_{j}=\frac{1}{\alpha}\left[\frac{a-1}{a-b}(X.Y)^{a-2}-\frac{b-1}{a-b}(X.Y)^{b-2}\right]\times\left(\frac{\sum_{i}p^{\alpha}_{i}q^{1-\alpha}_{i}}{\sum_{i}q_{i}}\right)^{\frac{1}{\alpha}-1}p^{\alpha}_{j}q^{-\alpha}_{j}\\&
V_{j}=\frac{1}{\alpha}\left[\frac{a-1}{a-b}(X.Y)^{a-2}-\frac{b-1}{a-b}(X.Y)^{b-2}\right]\times\left(\frac{\sum_{i}p^{\alpha}_{i}q^{1-\alpha}_{i}}{\sum_{i}q_{i}}\right)^{\frac{1}{\alpha}}
\end{split}
\end{equation}
\textbf{When $\alpha<0$, we will have:}
\begin{equation}
\begin{split}
&U_{j}=-\frac{1}{\alpha}\left[\frac{a-1}{a-b}(X.Y)^{a-2}-\frac{b-1}{a-b}(X.Y)^{b-2}\right]\times\left(\frac{\sum_{i}p^{\alpha}_{i}q^{1-\alpha}_{i}}{\sum_{i}q_{i}}\right)^{\frac{1}{\alpha}}\\&
V_{j}=-\frac{1}{\alpha}\left[\frac{a-1}{a-b}(X.Y)^{a-2}-\frac{b-1}{a-b}(X.Y)^{b-2}\right]\times\left(\frac{\sum_{i}p^{\alpha}_{i}q^{1-\alpha}_{i}}{\sum_{i}q_{i}}\right)^{\frac{1}{\alpha}-1}p^{\alpha}_{j}q^{-\alpha}_{j}
\end{split}
\end{equation}
We can easily verify that, as for all invariant divergences (logarithmic or not), we have: 
\begin{equation}
\sum_{j}q_{j}\frac{\partial L_{d}D_{\alpha}I(p\|q)}{\partial q_{j}}=0	
\end{equation}
For comparison, we remind here the results corresponding to the application of the natural logarithm on the invariant " $\alpha$ divergence " (\ref{alphadivinv}); these results have been mentioned in a previous work mentioned above.\\
Given the property (\ref{lognprod}), it follows:
\begin{equation}
LD_{\alpha}I(p\|q)=T\left(\log A-\log X-\log Y\right)	
\end{equation}
Then:
\begin{equation}
-\frac{\partial LD_{\alpha}I(p\|q)}{\partial q_{j}}=-T\left(\frac{1}{A}\frac{\partial A}{\partial q_{j}}-\frac{1}{X}\frac{\partial X}{\partial q_{j}}-\frac{1}{Y}\frac{\partial Y}{\partial q_{j}}\right)	
\end{equation}
Taking into account (\ref{derivpartiellesa}), we have:
\begin{equation}
-\frac{\partial LD_{\alpha}I(p\|q)}{\partial q_{j}}=-\frac{1}{\alpha}\left[-\frac{1}{X}\left(\sum_{i}p^{\frac{1}{\alpha}}_{i}q^{1-\alpha}_{i}\right)^{\frac{1}{\alpha}-1}p^{\alpha}_{j}q^{-\alpha}_{j}+\frac{1}{Y}\left(\sum_{i}q_{i}\right)^{-\frac{1}{\alpha}}\right]
\end{equation}
And finally, by expliciting "$X$" and "$Y$":
\begin{equation}
-\frac{\partial LD_{\alpha}I(p\|q)}{\partial q_{j}}=-\frac{1}{\alpha}\left[\left(\sum_{i}q_{i}\right)^{-1}-\left(\sum_{i}p^{\alpha}_{i}q^{1-\alpha}_{i}\right)^{-1}p^{\alpha}_{j}q^{-\alpha}_{j}\right]	
\end{equation}
This expression can also be obtained from the relation (\ref{gradalphadivinv}) by following the method detailed in Appendix 1.\\
The decomposition of this expression in the form $U-V$ changes according to the sign of "$\alpha$".\\
As for all the scale invariants divergences, we will have:
\begin{equation}
\sum_{j}q_{j}\frac{\partial LD_{\alpha}I(p\|q)}{\partial q_{j}}=0	
\end{equation}

\subsection{Beta ($\beta$) divergence.}
This divergence already analyzed in \cite{lanteri2019}, \cite{lanteri2020} and in the works cited there \cite{amari2009},\cite{cichocki2010},\cite{cichocki2011} is expressed in the form:
\begin{equation}
D_{\beta}(p\|q)=\frac{1}{\beta(\beta-1)}\left\{\sum_{i}p^{\beta}_{i}-\left[\sum_{i}\beta\;p_{i}q^{\beta-1}_{i}+(1-\beta)\;q^{\beta}_{i}\right]\right\}	
\end{equation}
It can be expressed in the equivalent form:
\begin{equation}
D_{\beta}(p\|q)=\frac{1}{\beta(\beta-1)}\sum_{i}p^{\beta}_{i}-\frac{1}{\beta-1}\sum_{i}p_{i}q^{\beta-1}_{i}+\frac{1}{\beta}q^{\beta}_{i}	
\end{equation}
As already mentioned in \cite{fevotte2011}, this divergence is convex for $1<\beta\leq2$, and we can further add that it is pseudo-convex in the other cases.\\
The decomposition of this divergence in the form $A-B$ with $A>0$ and $B>0$ implies to distinguish 3 cases, according to the values of the parameter "$\beta$".\\

\textbf{*Case 1:} $0<\beta<1$:
\begin{equation}
A=\frac{1}{1-\beta}\sum_{i}p_{i}q^{\beta-1}_{i}+\frac{1}{\beta}q^{\beta}_{i}\;\;;\;\;B=\frac{1}{\beta(1-\beta)}\sum_{i}p_{i}^{\beta}
\end{equation}
Taking into account the separability of the terms $A(p,q)$ and $B(p,q)$, it comes:
\begin{equation}
\begin{split}
&\frac{\partial A}{\partial q_{j}}=\frac{\partial A_{j}}{\partial q_{j}}=-p_{j}q^{\beta-2}_{j}+q^{\beta-1}_{j}
\\&\frac{\partial B}{\partial q_{j}}=\frac{\partial B_{j}}{\partial q_{j}}=0
\end{split}
\end{equation}
By application of the deformed logarithm on each term of the difference, we obtain (\ref{logdDexplicite}) whose gradient with respect to "$q$" is given by (\ref{gradlogdD}).\\
Hence the expression of the opposite of the gradient with respect to "$q$":
\begin{equation}
-\frac{\partial L_{d}D_{\beta}(p\|q)}{\partial q_{j}}=\left[\frac{a-1}{a-b}A^{a-2}-\frac{b-1}{a-b}A^{b-2}\right]\left[p_{j}q^{\beta-2}_{j}-q^{\beta-1}_{j}\right] 
\end{equation}
The decomposition of the opposite of the gradient in the form $U-V$ with $U>0$ and $V>0$ is expressed:
\begin{equation}
\begin{split}
&U_{j}=\left[\frac{a-1}{a-b}A^{a-2}-\frac{b-1}{a-b}A^{b-2}\right]p_{j}q^{\beta-2}_{j}\\&
V_{j}=\left[\frac{a-1}{a-b}A^{a-2}-\frac{b-1}{a-b}A^{b-2}\right]q^{\beta-1}_{j}
\end{split}
\end{equation}

\textbf{*Case 2:} $\beta>1$:
\begin{equation}
A=\frac{1}{\beta(\beta-1)}\sum_{i}p_{i}^{\beta}+\frac{1}{\beta}q^{\beta}_{i}\;\;;\;\;B=\frac{1}{\beta-1}\sum_{i}p_{i}q^{\beta-1}_{i}	
\end{equation}
Taking into account the separability of the terms $A(p,q)$ and $B(p,q)$, it comes:
\begin{equation}
\begin{split}
&\frac{\partial A}{\partial q_{j}}=\frac{\partial A_{j}}{\partial q_{j}}=q^{\beta-1}_{j}
\\&\frac{\partial B}{\partial q_{j}}=\frac{\partial B_{j}}{\partial q_{j}}=p_{j}q^{\beta-2}_{j}
\end{split}
\end{equation}
By application of the deformed logarithm on each term of the difference, we obtain (\ref{logdDexplicite}) whose gradient with respect to "$q$" is given by (\ref{gradlogdD}).\\
Hence the expression of the opposite of the gradient with respect to "$q$":
\begin{equation}
-\frac{\partial L_{d}D_{\beta}(p\|q)}{\partial q_{j}}=-\left[\frac{a-1}{a-b}A^{a-2}-\frac{b-1}{a-b}A^{b-2}\right]q^{\beta-1}_{j}+\left[\frac{a-1}{a-b}B^{a-2}-\frac{b-1}{a-b}B^{b-2}\right]p_{j}q^{\beta-2}_{j}
\end{equation}
The decomposition of the opposite of the gradient in the form $U-V$ with $U>0$ and $V>0$ is expressed:
\begin{equation}
\begin{split}
&U_{j}=\left[\frac{a-1}{a-b}B^{a-2}-\frac{b-1}{a-b}B^{b-2}\right]p_{j}q^{\beta-2}_{j}\\&
V_{j}=\left[\frac{a-1}{a-b}A^{a-2}-\frac{b-1}{a-b}A^{b-2}\right]q^{\beta-1}_{j}
\end{split}
\end{equation}

\textbf{*Case 3:} $\beta<0$.\\
This is not a classical situation, but it must be considered.\\
Here we have:
\begin{equation}
A=\frac{1}{\beta(\beta-1)}\sum_{i}p_{i}^{\beta}-\frac{1}{\beta-1}\sum_{i}p_{i}q^{\beta-1}_{i}	\;\;;\;\;B=-\frac{1}{\beta}q^{\beta}_{i}
\end{equation}
Taking into account the separability of the terms $A(p,q)$ and $B(p,q)$, it comes:
\begin{equation}
\begin{split}
&\frac{\partial A}{\partial q_{j}}=\frac{\partial A_{j}}{\partial q_{j}}=-p_{j}q^{\beta-2}_{j}
\\&\frac{\partial B}{\partial q_{j}}=\frac{\partial B_{j}}{\partial q_{j}}=-q^{\beta-1}_{j}
\end{split}
\end{equation}
By application of the deformed logarithm on each term of the difference, we obtain (\ref{logdDexplicite}) whose gradient with respect to "$q$" is given by (\ref{gradlogdD}).\\
Hence the expression of the opposite of the gradient with respect to "$q$":
\begin{equation}
-\frac{\partial L_{d}D_{\beta}(p\|q)}{\partial q_{j}}=\left[\frac{a-1}{a-b}A^{a-2}-\frac{b-1}{a-b}A^{b-2}\right]p_{j}q^{\beta-2}_{j}-\left[\frac{a-1}{a-b}B^{a-2}-\frac{b-1}{a-b}B^{b-2}\right]q^{\beta-1}_{j}
\end{equation}
The decomposition of the opposite of the gradient in the form $U-V$ with $U>0$ and $V>0$ is expressed:
\begin{equation}
\begin{split}
&U_{j}=\left[\frac{a-1}{a-b}A^{a-2}-\frac{b-1}{a-b}A^{b-2}\right]p_{j}q^{\beta-2}_{j}\\&
V_{j}=\left[\frac{a-1}{a-b}B^{a-2}-\frac{b-1}{a-b}B^{b-2}\right]q^{\beta-1}_{j}
\end{split}
\end{equation}

\subsection{Scale invariant Beta ($\beta$) divergence.}
This divergence has been studied in a previous work already cited \cite{lanteri2019}, \cite{lanteri2020}.\\
By introducing its nominal invariance factor which is explicitly derivable and is expressed:
 \begin{equation}
 K_{0,\beta}=\frac{\sum_{i}p_{i}q^{\beta-1}_{i}}{\sum_{i}q^{\beta}_{i}}
 \end{equation}
We obtain the invariant form:
\begin{equation}
D_{\beta}I(p\|q)=\underbrace{\frac{1}{\beta(\beta-1)}}_{T}\left[\underbrace{\sum_{i}p^{\beta}_{i}}_{A}-\underbrace{\left(\sum_{i}p_{i}q^{\beta-1}_{i}\right)^{\beta}}_{X}\times\underbrace{\left(\sum_{i}q^{\beta}_{i}\right)^{1-\beta}}_{Y}\right]	
\end{equation}
Since the terms "$X$" and "$Y$" are not separable, the calculation of the partial derivatives leads to:
\begin{equation}
\begin{split}
&\frac{\partial A}{\partial q_{j}}=0
\\&\frac{\partial X}{\partial q_{j}}=\beta(\beta-1)\left(\sum_{i}p_{i}q^{\beta-1}_{i}\right)^{\beta-1}p_{j}q^{\beta-2}_{j}
\\&\frac{\partial Y}{\partial q_{j}}=\beta(1-\beta)\left(\sum_{i}q^{\beta}_{i}\right)^{-\beta}.\ q^{\beta-1}_{j}
\label{derivpartiellesb}	
\end{split}	
\end{equation}
By application of the deformed logarithm, following the development of section (3.3.3), we obtain the expression of the divergence (\ref{logdD2bisdev}) and that of  the gradient with respect to "$q$" (\ref{gradlogdD2bis}).\\
The opposite of the gradient then expressed with (\ref{derivpartiellesb}):
\begin{equation}
	-\frac{\partial L_{d}D_{\beta}I(p\|q)}{\partial q_{j}}=T\left\{\underbrace{\left[\frac{a-1}{a-b}(X.Y)^{a-2}-\frac{b-1}{a-b}(X.Y)^{b-2}\right]}_{Z>0}\frac{\partial (X.Y)}{\partial q_{j}}\right\}
\end{equation}
With:
\begin{equation}
\frac{\partial (X.Y)}{\partial q_{j}}=\beta(1-\beta)\left[\left(\frac{\sum_{i}p_{i}q^{\beta-1}_{i}}{\sum_{i}q^{\beta}_{i}}\right)^{\beta}q^{\beta-1}_{j}-\left(\frac{\sum_{i}p_{i}q^{\beta-1}_{i}}{\sum_{i}q^{\beta}_{i}}\right)^{\beta-1}p_{j}q^{\beta-2}_{j}\right]	
\end{equation}
And finally, taking into account the expression of "$T$", the opposite of the gradient is expressed:
\begin{equation}
-\frac{\partial L_{d}D_{\beta}I(p\|q)}{\partial q_{j}}=Z\;\left[\left(\frac{\sum_{i}p_{i}q^{\beta-1}_{i}}{\sum_{i}q^{\beta}_{i}}\right)^{\beta-1}p_{j}q^{\beta-2}_{j}-\left(\frac{\sum_{i}p_{i}q^{\beta-1}_{i}}{\sum_{i}q^{\beta}_{i}}\right)^{\beta}q^{\beta-1}_{j}\right]	
\end{equation}

The decomposition of the opposite of the gradient in the form $U-V$ with $U>0$ and $V>0$ is straightforward, and we will have as for all invariant divergences:
\begin{equation}
\sum_{j}q_{j}\frac{\partial L_{d}BEI(p\|q)}{\partial q_{j}}=0	
\end{equation}

\subsection{Alpha-Beta ($\alpha\beta$) divergence.}
This divergence already analyzed in \cite{lanteri2019}, \cite{lanteri2020} and in the papers cited there \cite{amari2009},\cite{cichocki2010},\cite{cichocki2011} is expressed:
\begin{equation}
D_{\alpha\beta}(p\|q)=T\left[\sum_{i}p^{\alpha+\beta-1}_{i}+\frac{\beta-1}{\alpha}\sum_{i}q^{\alpha+\beta-1}_{i}-\frac{\alpha+\beta-1}{\alpha}\sum_{i}p^{\alpha}_{i}q^{\beta-1}_{i}\right]
\label{divab}   
\end{equation}

With:
\begin{equation}
T=\frac{1}{(\beta-1)(\alpha+\beta-1)}
\end{equation}
We can write in an equivalent way:
\begin{equation}
\begin{split}
D_{\alpha\beta}(p\|q)=&\frac{1}{(\beta-1)(\alpha+\beta-1)}\sum_{i}p^{\alpha+\beta-1}_{i}\\&+\frac{1}{\alpha(\alpha+\beta-1)}\sum_{i}q^{\alpha+\beta-1}_{i}\\&-\frac{1}{\alpha(\beta-1)}\sum_{i}p^{\alpha}_{i}q^{\beta-1}_{i}
\label{divabbis}
\end{split}   
\end{equation}
Notice that the variable "$\beta$" of the reference \cite{cichocki2011} is replaced by $(\beta-1)$ in our expression.\\
Writing this divergence in the form $A-B$, $A>0$, $B>0$ and applying the deformed logarithm, we obtain:
\begin{equation}
L_{d}D_{\alpha\beta}(p\|q)=\frac{1}{a-b}\left[A^{a-1}-A^{b-1}-\left( B^{a-1}-B^{b-1}\right) \right]
\end{equation}
The general expression for the opposite of the gradient with respect to "$q$" is written, $\forall j$:
\begin{equation}
\begin{split}
-\frac{\partial L_{d}D_{\alpha\beta}(p\|q)}{\partial q_{j}}=&-\left[\left(\underbrace{\frac{a-1}{a-b}A^{a-2}-\frac{b-1}{a-b}A^{b-2}}_{Z_{A}>0}\right)\frac{\partial A}{\partial q_{j}}\right]\\&+\left[\left(\underbrace{\frac{a-1}{a-b}B^{a-2}-\frac{b-1}{a-b}B^{b-2}}_{Z_{B}>0}\right)\frac{\partial B}{\partial q_{j}}\right]
\end{split}  
\end{equation}
The decomposition of the divergence $D_{\alpha\beta}(p|q)$ in the form $(A-B)$ with $A$ and $B$ positive depends on the signs of $\alpha$, $\beta-1$ and $\alpha+\beta-1$; this implies a particular analysis.\\
\textbf{If $\alpha>0$, four situations can occur:}\\

For each of the situations considered, the decomposition of the opposite of the gradient in the form $U-V$ with $U>0$ and $V>0$ is indicated.\\

\textbf{*1 -} \;$\alpha>0$, $\beta-1>0$ then $\alpha+\beta-1>0$, we have:\\
\begin{equation}
\begin{split}
&A=\frac{1}{(\beta-1)(\alpha+\beta-1)}\sum_{i}p^{\alpha+\beta-1}_{i}+\frac{1}{\alpha(\alpha+\beta-1)}\sum_{i}q^{\alpha+\beta-1}_{i}\\&
B=\frac{1}{\alpha(\beta-1)}\sum_{i}p^{\alpha}_{i}q^{\beta-1}_{i}
\end{split}
\end{equation}
Then:
\begin{equation}
\begin{split}
&\frac{\partial A}{\partial q_{j}}=\frac{1}{\alpha}q^{\alpha+\beta-2}_{j}\\&
\frac{\partial B}{\partial q_{j}}=\frac{1}{\alpha}p^{\alpha}_{j}q^{\beta-2}_{j}
\end{split}
\end{equation}
And we will have:
\begin{equation}
-\frac{\partial L_{d}D_{\alpha\beta}(p\|q)}{\partial q_{j}}=\underbrace{\frac{Z_{B}}{\alpha}p^{\alpha}_{j}q^{\beta-2}_{j}}_{U_{j}}-\underbrace{\frac{Z_{A}}{\alpha}q^{\alpha+\beta-2}_{j}}_{V_{j}}
\end{equation}

\textbf{*2 -} \;$\alpha>0$, $\beta-1<0$ and $\alpha+\beta-1>0$, we then have:\\
\begin{equation}
\begin{split}
&A=\frac{1}{\alpha(\alpha+\beta-1)}\sum_{i}q^{\alpha+\beta-1}_{i}+\frac{1}{\alpha(1-\beta)}\sum_{i}p^{\alpha}_{i}q^{\beta-1}_{i}\\&
B=\frac{1}{(1-\beta)(\alpha+\beta-1)}\sum_{i}p^{\alpha+\beta-1}_{i}
\end{split}
\end{equation}
Then, we will have:
\begin{equation}
\begin{split}
&\frac{\partial A}{\partial q_{j}}=\frac{1}{\alpha}q^{\alpha+\beta-2}_{j}-\frac{1}{\alpha}p^{\alpha}_{j}q^{\beta-2}_{j}\\&
\frac{\partial B}{\partial q_{j}}=0
\end{split}
\end{equation}
Consequently:
\begin{equation}
-\frac{\partial L_{d}D_{\alpha\beta}(p\|q)}{\partial q_{j}}=\underbrace{\frac{Z_{A}}{\alpha}p^{\alpha}_{j}q^{\beta-2}_{j}}_{U_{j}}-\underbrace{\frac{Z_{A}}{\alpha}q^{\alpha+\beta-2}_{j}}_{V_{j}}
\end{equation}

\textbf{*3 -} \;$\alpha>0$, $\beta-1>0$ and $\alpha+\beta-1<0$, this case is not to be considered, because it is obviously impossible.\\

\textbf{*4 -} \;$\alpha>0$, $\beta-1<0$ and $\alpha+\beta-1<0$, we have:\\
\begin{equation}
\begin{split}
&A=\frac{1}{(\beta-1)(\alpha+\beta-1)}\sum_{i}p^{\alpha+\beta-1}_{i}-\frac{1}{\alpha(\beta-1)}\sum_{i}p^{\alpha}_{i}q^{\beta-1}_{i}\\&
B=\frac{1}{\alpha(1-\alpha-\beta)}\sum_{i}q^{\alpha+\beta-1}_{i}
\end{split}
\end{equation}
This gives us:
\begin{equation}
\begin{split}
&\frac{\partial A}{\partial q_{j}}=-\frac{1}{\alpha}p^{\alpha}_{j}q^{\beta-2}_{j}\\&
\frac{\partial B}{\partial q_{j}}=-\frac{1}{\alpha}q^{\alpha+\beta-2}_{j}
\end{split}
\end{equation}
Then:
\begin{equation}
-\frac{\partial L_{d}D_{\alpha\beta}(p\|q)}{\partial q_{j}}=\underbrace{\frac{Z_{A}}{\alpha}p^{\alpha}_{j}q^{\beta-2}_{j}}_{U_{j}}-\underbrace{\frac{Z_{B}}{\alpha}q^{\alpha+\beta-2}_{j}}_{V_{j}}
\end{equation}
\textbf{If, now, we consider the cases corresponding to $\alpha<0$, the 4 previous situations must be examined again.}\\

\textbf{*1bis -} \;$\alpha<0$, $\beta-1>0$ and $\alpha+\beta-1>0$; the expressions de $A$, $B$, $\frac{\partial A}{\partial q_{j}}$, $\frac{\partial B}{\partial q_{j}}$  are analogous to the \textbf{(case *4)} above, however, taking into account the sign of "$\alpha$", the decomposition of the opposite of the gradient of the divergence is modified and is written in the form $U-V$, $U>0$, $V>0$ by taking:
\begin{equation}
U_{j}=-\frac{Z_{B}}{\alpha}q^{\alpha+\beta-2}_{j}\;\;;\;\; V_{j}=-\frac{Z_{A}}{\alpha}p^{\alpha}_{j}q^{\beta-2}_{j}
\end{equation}
\textbf{*2bis -} \;$\alpha<0$, $\beta-1<0$ and $\alpha+\beta-1>0$; this case is clearly impossible.\\\\
\textbf{*3bis -} \;$\alpha<0$, $\beta-1>0$ and $\alpha+\beta-1<0$; the expressions of $A$, $B$, $\frac{\partial A}{\partial q_{j}}$, $\frac{\partial B}{\partial q_{j}}$ are analogous to the \textbf{(case *2)} above, however, given the sign of "$\alpha$", the decomposition of the opposite of the gradient of the divergence is modified and can be written as $U-V$, $U>0$, $V>0$ by taking:
\begin{equation}
U_{j}=-\frac{Z_{A}}{\alpha}q^{\alpha+\beta-2}_{j}\;\;;\;\; V_{j}=-\frac{Z_{A}}{\alpha}q^{\alpha+\beta-2}_{j}
\end{equation}

\textbf{*4bis -} \;$\alpha<0$, $\beta-1<0$ then $\alpha+\beta-1<0$; the expressions of $A$, $B$, $\frac{\partial A}{\partial q_{j}}$, $\frac{\partial B}{\partial q_{j}}$ are analogous to the \textbf{(case *1)} above, however, given the sign of "$\alpha$", the decomposition of the opposite of the gradient of the divergence is modified and can be written as $U-V$, $U>0$, $V>0$ by taking:
\begin{equation}
U_{j}=-\frac{Z_{A}}{\alpha}q^{\alpha+\beta-2}_{j}\;\;;\;\; V_{j}=-\frac{Z_{B}}{\alpha}p^{\alpha}_{j}q^{\beta-2}_{j}
\end{equation}

\subsection{Scale invariant Alpha-Beta ($\alpha\beta$) divergence.}
Given the $(\alpha\beta)$ divergence (\ref{divabbis}),
the nominal invariance factor can be calculated explicitly; it is written:
\begin{equation}
K_{0,\alpha\beta}=\left[\frac{\sum_{i}p^{\alpha}_{i}q^{\beta-1}_{i}}{\sum_{i}q^{\alpha+\beta-1}_{i}} \right]^{\frac{1}{\alpha}} 
\end{equation}
Introducing this invariance factor in (\ref{divabbis}), we obtain the ($\alpha\beta$) invariant divergence already considered in a previous work \cite{lanteri2019}, \cite{lanteri2020}, which is written:
\begin{equation}
D_{\alpha\beta}I(p\|q)=T\;\left[\underbrace{\sum_{i}p^{\alpha+\beta-1}_{i}}_{A(p,q)}-\underbrace{\left(\sum_{i}p^{\alpha}_{i}q^{\beta-1}_{i}\right)^{\frac{\alpha+\beta-1}{\alpha}}}_{X(p,q)}\times\underbrace{\left(\sum_{i}q^{\alpha+\beta-1}_{i}\right)^{\frac{1-\beta}{\alpha}}}_{Y(p,q)}\right]	
\end{equation}
With:
\begin{equation}
T=\frac{1}{(\beta-1)(\alpha+\beta-1)}	
\end{equation}
By applying the deformed Logarithm, it comes:
\begin{equation}
L_{d}D_{\alpha\beta}I(p\|q)=T\;\left\{\frac{A^{a-1}-A^{b-1}}{a-b}-\left[\frac{(X.Y)^{a-1}-(X.Y)^{b-1}}{a-b}\right]\right\}	
\end{equation}
With:
\begin{equation}
\begin{split}
&\frac{\partial A}{\partial q_{j}}=\frac{\partial A_{j}}{\partial q_{j}}=0\\&
\frac{\partial X}{\partial q_{j}}=\frac{(\alpha+\beta-1)(\beta-1)}{\alpha}\left(\sum_{i}p^{\alpha}_{i}q^{\beta-1}_{i}\right)^{\frac{\beta-1}{\alpha}}p^{\alpha}_{j}q^{\beta-2}_{j}\\&
\frac{\partial Y}{\partial q_{j}}=\frac{(\alpha+\beta-1)(1-\beta)}{\alpha}\left(\sum_{i}q^{\alpha+\beta-1}_{i}\right)^{\frac{1-\alpha-\beta}{\alpha}}q^{\alpha+\beta-2}_{j}
\end{split}
\end{equation}
The opposite of the gradient with respect to "$q$" is written $\forall j$:
\begin{equation}
	-\frac{\partial L_{d}D_{\alpha\beta}I(p\|q)}{\partial q_{j}}=T\;\left[\left(\frac{a-1}{a-b}\right)(X.Y)^{a-2}-\left(\frac{b-1}{a-b}\right)(X.Y)^{b-2}\right]\frac{\partial (X.Y)}{\partial q_{j}}
\end{equation}
With:
\begin{equation}
\begin{split}
\frac{\partial (X.Y)}{\partial q_{j}}=&\frac{(1-\beta)(\alpha+\beta-1)}{\alpha}\underbrace{\left(\frac{\sum_{i}p^{\alpha}_{i}q^{\beta-1}_{i}}{\sum_{i}q^{\alpha+\beta-1}_{i}}\right)^{\frac{\alpha+\beta-1}{\alpha}}q^{\alpha+\beta-2}_{j}}_{M>0}\\&
-\frac{(1-\beta)(\alpha+\beta-1)}{\alpha}\underbrace{\left(\frac{\sum_{i}p^{\alpha}_{i}q^{\beta-1}_{i}}{\sum_{i}q^{\alpha+\beta-1}_{i}}\right)^{\frac{\beta-1}{\alpha}}p^{\alpha}_{j}q^{\beta-2}_{j}}_{N>0}
\end{split}                          
\end{equation}
Hence, finally, taking into account the expression of "$T$", it follows:
\begin{equation}
	-\frac{\partial L_{d}D_{\alpha\beta}I(p\|q)}{\partial q_{j}}=-\frac{1}{\alpha}\left[\left(\frac{a-1}{a-b}\right)(X.Y)^{a-2}-\left(\frac{b-1}{a-b}\right)(X.Y)^{b-2}\right]\left[M-N\right]
\end{equation}
We can observe that this expression does not depend on the relative values of "$\alpha$" and "$\beta$" which makes any discussion on this point useless.\\
A discussion according to the sign of "$\alpha$" must however be carried out, indeed,
the decomposition of the opposite of the gradient in the form $U-V$ with $U>0$ and $V>0$ will depend on the sign of "$\alpha$", and we will have:\\
If $\alpha>0$:
\begin{equation}
\begin{split}
&U_{j}=\frac{1}{\alpha}\left[\left(\frac{a-1}{a-b}\right)(X.Y)^{a-2}-\left(\frac{b-1}{a-b}\right)(X.Y)^{b-2}\right]\times N\\&V_{j}=\frac{1}{\alpha}\left[\left(\frac{a-1}{a-b}\right)(X.Y)^{a-2}-\left(\frac{b-1}{a-b}\right)(X.Y)^{b-2}\right]\times M
\end{split}                          
\end{equation}
If $\alpha<0$:
\begin{equation}
\begin{split}
&U_{j}=-\frac{1}{\alpha}\left[\left(\frac{a-1}{a-b}\right)(X.Y)^{a-2}-\left(\frac{b-1}{a-b}\right)(X.Y)^{b-2}\right]\times M\\&V_{j}=-\frac{1}{\alpha}\left[\left(\frac{a-1}{a-b}\right)(X.Y)^{a-2}-\left(\frac{b-1}{a-b}\right)(X.Y)^{b-2}\right]\times N
\end{split}                          
\end{equation}

\section{Application to divergences between means.}
These divergences have been developed in Chapter 7 of the references \cite{lanteri2019}, \cite{lanteri2020}, they are based on the previous works of Taneja \cite{taneja2001} and Ben-Tal \cite{ben1989}; they have all the form "$D$" (\ref{divD}); we analyze here only some specific cases as examples.\\
It should be noted that a generalization of these divergences has been proposed by Taneja in a form using the Tsallis deformed Logarithm; we will remind its main points in Appendix 2.

\subsection{Geometric-Harmonic Divergence.}
With $0\leq \alpha\leq 1$:
\begin{equation}
\begin{split}
&MG=\sum_{i} p^{\alpha}_{i}q^{1-\alpha}_{i}=\sum_{i}MG_{i}\\&
MH=\sum_{i}\frac{p_{i}q_{i}}{\left(1-\alpha\right)p_{i}+\alpha q_{i}}=\sum_{i}MH_{i}
\end{split}
\end{equation}
The divergence is written:
\begin{equation}
DGH(p\|q)=MG-MH
\end{equation}
By applying the deformed logarithm, we obtain:
\begin{equation}
L_{d}DGH(p\|q)=\frac{\left(MG\right)^{a-1}-\left(MG\right)^{b-1}}{a-b}-\left[\frac{\left(MH\right)^{a-1}-\left(MH\right)^{b-1}}{a-b}\right] 
\end{equation}
The opposite of the gradient with respect to "$q$" is written $\forall j$:
\begin{equation}
\begin{split}
-\frac{\partial L_{d}DGH(p\|q)}{\partial q_{j}}=&-\left[\underbrace{\frac{a-1}{a-b}\left(MG\right)^{a-2}-\frac{b-1}{a-b}\left(MG\right)^{b-2}}_{Z}\right]\frac{\partial \left(MG\right)}{\partial q_{j}}\\&+\left[\underbrace{\frac{a-1}{a-b}\left(MH\right)^{a-2}-\frac{b-1}{a-b}\left(MH\right)^{b-2}}_{Y}\right]\frac{\partial \left(MH\right)}{\partial q_{j}}
\end{split}
\end{equation}
With:
\begin{equation}
\begin{split}
&\frac{\partial \left(MG\right)}{\partial q_{j}}=\frac{\partial \left(MG_{j}\right)}{\partial q_{j}}=\left(1-\alpha\right)\frac{p^{\alpha}_{j}}{q^{\alpha}_{j}}>0\\&\frac{\partial \left(MH\right)}{\partial q_{j}}=\frac{\partial \left(MH_{j}\right)}{\partial q_{j}}=\frac{\left(1-\alpha\right)p^{2}_{j}}{\left[\left(1-\alpha\right)p_{j}+\alpha q_{j}\right]^{2}}>0
\end{split}
\end{equation}
The decomposition of the opposite of the gradient in the form $U-V$ with $U>0$ and $V>0$ is straightforward.

\subsection{Scale invariant Geometric-Harmonic Divergence.}
For this divergence, the nominal invariance factor has no explicit expression, so we will use the invariance factor $K^{*}=\frac{sum_{j}p_{j}}{\sum_{j}q_{j}}$.\\
By using the normalized variables $\overline{p}_{i}=\frac{p_{i}}{\sum_{j}p_{j}}$ and $\overline{q}_{i}=\frac{q_{i}}{\sum_{j}q_{j}}$, the divergence in invariant form is written:
\begin{equation}
DGHI(p\|q)=\sum_{j}p_{j}\left[\overline{MG}-\overline{MH}\right]  
\end{equation}
With:
\begin{equation}
\begin{split}
&\overline{MG}=\sum_{i}\overline{MG}_{i}=\sum_{i}\overline{p}^{\alpha}_{i}\overline{q}^{1-\alpha}_{i}
\\&\overline{MH}=\sum_{i}\overline{MH}_{i}=\sum_{i}\frac{\overline{p}_{i}\overline{q}_{i}}{\left(1-\alpha\right)\overline{p}_{i}+\alpha\overline{q}_{i}}
\end{split}
\end{equation}
By applying the deformed logarithm, we obtain:
\begin{equation}
\begin{split}
L_{d}DGHI(p\|q)=&\frac{\sum_{j}p_{j}}{a-b}\left[\left(\overline{MG}\right)^{a-1}-\left(\overline{MG}\right)^{b-1}\right]
\\&-\frac{\sum_{j}p_{j}}{a-b}\left[\left(\overline{MH}\right)^{a-1}-\left(\overline{MH}\right)^{b-1}\right] 
\end{split}   
\end{equation}
The opposite of the gradient with respect to "$q$" is written $\forall j$:
\begin{equation}
\begin{split}
-\frac{\partial L_{d}DGHI(p\|q)}{\partial q_{j}}&=-\sum_{j}p_{j}\left[\underbrace{\frac{a-1}{a-b}\left(\overline{MG}\right)^{a-2}-\frac{b-1}{a-b}\left(\overline{MG}\right)^{b-2}}_{Z}\right]\frac{\partial \overline{MG}}{\partial q_{j}}\\&
+\sum_{j}p_{j}\left[\underbrace{\frac{a-1}{a-b}\left(\overline{MH}\right)^{a-2}-\frac{b-1}{a-b}\left(\overline{MH}\right)^{b-2}}_{Y}\right]\frac{\partial \overline{MH}}{\partial q_{j}}
\end{split}
\end{equation}
With:
\begin{equation}
\frac{\partial \overline{MG}}{\partial q_{j}}=\frac{1-\alpha}{\sum_{j}q_{j}}\left[\overline{p}^{\alpha}_{j}\overline{q}^{-\alpha}_{j}-\sum_{i}\overline{p}^{\alpha}_{i}\overline{q}^{1-\alpha}_{i}\right] =\frac{1-\alpha}{\sum_{j}q_{j}}\left[ \frac{\overline{MG}_{j}}{\overline{q}_{j}}-\sum_{i}\overline{MG}_{i}\right] 
\end{equation}
And:
\begin{equation}
\begin{split}
\frac{\partial \overline{MH}}{\partial q_{j}}=&\frac{1-\alpha}{\sum_{j}q_{j}}\left\lbrace \frac{\overline{p}^{2}_{j}}{\left[\left(1-\alpha\right)\overline{p}_{j}+\alpha\overline{q}_{j}\right]^{2}}
-\sum_{i}\overline{q}_{i} \frac{ \overline{p}^{2}_{i}}{\left[ \left(1-\alpha\right)\overline{p}_{i}+\alpha\overline{q}_{i}\right]^{2}}\right\rbrace  \\&=\frac{1-\alpha}{\sum_{j}q_{j}}\left[ \frac{\overline{MH}^{2}_{j}}{\overline{q}^{2}_{j}}-\sum_{i}\frac{\overline{MH}^{2}_{i}}{\overline{q}_{i}}\right] 
\end{split}
\end{equation}
The decomposition of the opposite of the gradient in the form $U-V$ with $U>0$ and $V>0$ is:
\begin{equation}
\begin{split}
&U_{j}=(1-\alpha)\frac{\sum_{j}p_{j}}{\sum_{j}q_{j}}\left[ Y\frac{\overline{MH}^{2}_{j}}{\overline{q}^{2}_{j}}+Z\sum_{i}\overline{MG}_{i}\right]\\&
V_{j}=(1-\alpha)\frac{\sum_{j}p_{j}}{\sum_{j}q_{j}}\left[ Z\frac{\overline{MG}_{j}}{\overline{q}_{j}}+Y\sum_{i}\frac{\overline{MH}^{2}_{i}}{\overline{q}_{i}}\right]
\end{split}
\end{equation}
As for all invariant divergence, we will have:
\begin{equation}
\sum_{j}q_{j}\frac{\partial L_{d}DGHI(p\|q)}{\partial q_{j}}=0
\end{equation}

\subsection{Arithmetic-Geometric Divergence.}
Except for a multiplicative factor, this divergence is identical to the Alpha divergence discussed in section (6.1); consequently, we will not return to this point.\\

\subsection{Scale invariant Arithmetic-Geometric Divergence.}
In this case, although the nominal invariance factor is explicitly computable and leads to the invariant Alpha divergence already analyzed in section (6.2), we analyze for the sake of homogeneity, another form of invariant Arithmetic-Geometric divergence obtained by using the invariance factor $K^{*}=\frac{\sum_{j}p_{j}}{\sum_{j}q_{j}}$.\\
With:
\begin{equation}
\overline{MG}=\sum_{i} \overline{p}^{\alpha}_{i}\overline{q}^{1-\alpha}_{i}\:\:;\:\overline{MA}=\sum_{i}\alpha \overline{p}_{i}+\left(1-\alpha\right)\overline{q}_{i}=1
\end{equation}
The invariant divergence is written:
\begin{equation}
DAGI(p\|q)=\sum_{j}p_{j}\left[\overline{MA}-\overline{MG}\right]=\sum_{j}p_{j}\left[1-\overline{MG}\right] 
\end{equation}
By applying the deformed logarithm, we obtain:
\begin{equation}
L_{d}DAGI(p\|q)=-\sum_{j}p_{j}\left[\frac{\left(\overline{MG}\right)^{a-1}-\left(\overline{MG}\right)^{b-1}}{a-b}\right] 
\end{equation}
The opposite of the gradient with respect to "$q$" is written $\forall j$:
\begin{equation}
-\frac{\partial L_{d}DAGI(p\|q)}{\partial q_{j}}=\sum_{j}p_{j}\left[\underbrace{\frac{a-1}{a-b}\left(\overline{MG}\right)^{a-2}-\frac{b-1}{a-b}\left(\overline{MG}\right)^{b-2}}_{Z}\right]\frac{\partial \left(\overline{MG}\right)}{\partial q_{j}}
\end{equation}
With:
\begin{equation}
\frac{\partial \overline{MG}}{\partial q_{j}}=\frac{\left(1-\alpha\right)}{\sum_{j}q_{j}}\left[\overline{p}^{\alpha}_{j}\overline{q}^{-\alpha}_{j}-\sum_{i}\overline{p}^{\alpha}_{i}\overline{q}^{1-\alpha}_{i}\right]=\frac{\left(1-\alpha\right)}{\sum_{j}q_{j}}\left[\frac{\overline{MG}_{j}}{\overline{q}_{j}}-\overline{MG}\right]
\end{equation}
The decomposition of the opposite of the gradient in the form $U-V$, with $U>0$ and $V>0$ is written:
\begin{equation}
\begin{split}
&U_{j}=(1-\alpha)\frac{\sum_{j}p_{j}}{\sum_{j}q_{j}}Z\;\frac{\overline{MG}_{j}}{\overline{q}_{j}}\\&
V_{j}=(1-\alpha)\frac{\sum_{j}p_{j}}{\sum_{j}q_{j}}Z\;\overline{MG}
\end{split}
\end{equation}
As for all scale invariant divergence, we will have:
\begin{equation}
\sum_{j}q_{j}\frac{\partial L_{d}DAGI(p\|q)}{\partial q_{j}}=0
\end{equation}

\subsection{Arithmetic-Harmonic Divergence.}
With:
\begin{equation}
\begin{split}
&MA=\sum_{i} \alpha p_{i}+\left(1-\alpha\right)q_{i}=\sum_{i}MA_{i}\\&
MH=\sum_{i}\frac{p_{i}q_{i}}{\left(1-\alpha\right)p_{i}+\alpha q_{i}}=\sum_{i}MH_{i}
\end{split}
\end{equation}
The divergence is written:
\begin{equation}
DAH(p\|q)=MA-MH
\end{equation}
By applying the deformed logarithm, we obtain:
\begin{equation}
L_{d}DAH(p\|q)=\frac{\left(MA\right)^{a-1}-\left(MA\right)^{b-1}}{a-b}-\left[\frac{\left(MH\right)^{a-1}-\left(MH\right)^{b-1}}{a-b}\right] 
\end{equation}
The opposite of the gradient with respect to "$q$" is written $\forall j$:
\begin{equation}
\begin{split}
-\frac{\partial L_{d}DAH(p\|q)}{\partial q_{j}}=&-\left[\underbrace{\frac{a-1}{a-b}\left(MA\right)^{a-2}-\frac{b-1}{a-b}\left(MA\right)^{b-2}}_{Z>0}\right]\frac{\partial \left(MA\right)}{\partial q_{j}}\\&+\left[\underbrace{\frac{a-1}{a-b}\left(MH\right)^{a-2}-\frac{b-1}{a-b}\left(MH\right)^{b-2}}_{{Y>0}}\right]\frac{\partial \left(MH\right)}{\partial q_{j}}
\end{split}
\end{equation}
With:
\begin{equation}
\begin{split}
&\frac{\partial \left(MA\right)}{\partial q_{j}}=\frac{\partial \left(MA_{j}\right)}{\partial q_{j}}=\left(1-\alpha\right)>0\\&\frac{\partial \left(MH\right)}{\partial q_{j}}=\frac{\partial \left(MH_{j}\right)}{\partial q_{j}}=\frac{\left(1-\alpha\right)p^{2}_{j}}{\left[\left(1-\alpha\right)p_{j}+\alpha q_{j}\right]^{2}}>0
\end{split}
\end{equation}
The decomposition of the opposite of the gradient in the form $U-V$, with $U>0$ and $V>0$ is written immediately:
\begin{equation}
\begin{split}
&U_{j}=\frac{\left(1-\alpha\right)p^{2}_{j}}{\left[\left(1-\alpha\right)p_{j}+\alpha q_{j}\right]^{2}}\;Y\\&V_{j}=\left(1-\alpha\right) Z
\end{split}
\end{equation}

\subsection{Scale invariant Arithmetic-Harmonic Divergence.}
For this divergence, the nominal invariance factor has no explicit expression, so we will use the invariance factor $K^{*}=\frac{\sum_{j}p_{j}}{\sum_{j}q_{j}}$.\\
Using the normalized variables $\overline{p}_{i}=\frac{p_{i}}{\sum_{j}p_{j}}$ and $\overline{q}_{i}=\frac{q_{i}}{\sum_{j}q_{j}}$, we have:
\begin{equation}
\overline{MH}=\sum_{i} \frac{\overline{p}_{i}\overline{q}_{i}}{\left(1-\alpha\right)\overline{p}_{i}+\alpha \overline{q}_{i}}\:\:;\:\overline{MA}=\sum_{i}\alpha \overline{p}_{i}+\left(1-\alpha\right)\overline{q}_{i}=1
\end{equation}
The divergence is written:
\begin{equation}
DAHI(p\|q)=\sum_{j}p_{j}\left(\overline{MA}-\overline{MH}\right)=\sum_{j}p_{j}\left(1-\overline{MH}\right) 
\end{equation}
By applying the deformed logarithm, we obtain:
\begin{equation}
L_{d}DAHI(p\|q)=-\sum_{j}p_{j}\left[\frac{\left(\overline{MH}\right)^{a-1}-\left(\overline{MH}\right)^{b-1}}{a-b}\right] 
\end{equation}
The opposite of the gradient with respect to "$q$" is written $\forall j$:
\begin{equation}
-\frac{\partial L_{d}DAHI(p\|q)}{\partial q_{j}}=
\sum_{j}p_{j}\left[\underbrace{\frac{a-1}{a-b}\left(\overline{MH}\right)^{a-2}-\frac{b-1}{a-b}\left(\overline{MH}\right)^{b-2}}_{Z>0}\right]\frac{\partial \overline{MH}}{\partial q_{j}}
\end{equation}
With :
\begin{equation}
\begin{split}
\frac{\partial \overline{MH}}{\partial q_{j}}=&\frac{1-\alpha}{\sum_{j}q_{j}}\left\lbrace \frac{\overline{p}^{2}_{j}}{\left[\left(1-\alpha\right)\overline{p}_{j}+\alpha\overline{q}_{j}\right]^{2}}
-\sum_{i}\overline{q}_{i}\left[ \frac{\overline{p}_{i}}{\left(1-\alpha\right)\overline{p}_{i}+\alpha\overline{q}_{i}}\right]^{2}\right\rbrace  \\&=\frac{1-\alpha}{\sum_{j}q_{j}}\left[ \frac{\overline{MH}^{2}_{j}}{\overline{q}^{2}_{j}}-\sum_{i}\frac{\overline{MH}^{2}_{i}}{\overline{q}_{i}}\right] 
\end{split}
\end{equation}
As for all scale invariant divergence, we will have:
\begin{equation}
\sum_{j}q_{j}\frac{\partial L_{d}DAHI(p\|q)}{\partial q_{j}}=0
\end{equation}
The decomposition of the opposite of the gradient in the form $U-V$, with $U>0$ and $V>0$ is written immediately:
\begin{equation}
\begin{split}
&U_{j}=(1-\alpha)\frac{\sum_{j}p_{j}}{\sum_{j}q_{j}} Z \frac{\overline{MH}^{2}_{j}}{\overline{q}^{2}_{j}} \\&V_{j}=(1-\alpha)\frac{\sum_{j}p_{j}}{\sum_{j}q_{j}} Z \sum_{i}\frac{\overline{MH}^{2}_{i}}{\overline{q}_{i}} 
\end{split}
\end{equation}

\section{Application to F and G divergences.}
These divergences have the particularity of being based on the use of the Kullback-Leibler Divergence; their expressions therefore explicitly contain the natural Logarithm function.\\
The introduction of the deformed logarithm consists therefore simply in substituting in these expressions the natural logarithm by the deformed logarithm.

\subsection{F Divergence.}

It is written \cite{taneja2001}, \cite{lanteri2019}, \cite{lanteri2020}:
\begin{equation}
F\left(p\|q\right)=\sum_{i}p_{i}\log\frac{p_{i}}{\alpha p_{i}+\left(1-\alpha\right)q_{i}}+\left(1-\alpha\right)\left(q_{i}-p_{i}\right) 
\end{equation}
With the notation:
\begin{equation}
Z_{i}=\frac{p_{i}}{\alpha p_{i}+\left(1-\alpha\right)q_{i}}
\end{equation}
By introducing the deformed logarithm, it comes:
\begin{equation}
L_{d}F\left(p\|q\right)=\sum_{i}\frac{p_{i}}{a-b}\left[\left(Z_{i}\right)^{a-1}-\left(Z_{i}\right)^{b-1}\right]+\left(1-\alpha\right)\left(q_{i}-p_{i}\right)
\end{equation}
\begin{equation}
\frac{\partial L_{d}F\left(p\|q\right)}{\partial q_{j}}=\left[\frac{a-1}{a-b}p_{j}\left(Z_{j}\right)^{a-2}-\frac{b-1}{a-b}p_{j}\left(Z_{j}\right)^{b-2}\right]\frac{\partial Z_{j}}{\partial q_{j}}+\left(1-\alpha\right) 
\end{equation}
With:
\begin{equation}
\frac{\partial Z_{j}}{\partial q_{j}}=-\left(1-\alpha\right)\frac{p_{j}}{\left[\alpha p_{j}+\left(1-\alpha\right)q_{j}\right]^{2}}
\end{equation}
And finally, the expression of the opposite of the gradient with respect to "$q$" is written $\forall j$:
\begin{equation}
-\frac{\partial L_{d}F\left(p\|q\right)}{\partial q_{j}}=\left(1-\alpha\right)\left\lbrace\left[\frac{a-1}{a-b}\left(Z_{j}\right)^{a}-\frac{b-1}{a-b}\left(Z_{j}\right)^{b}\right]-1\right\rbrace 
\end{equation}
The decomposition of the opposite of the gradient in the form $U-V$, with $U>0$ and $V>0$ is written:
\begin{equation}
\begin{split}
&U_{j}=\left(1-\alpha\right)\left[\frac{a-1}{a-b}\left(Z_{j}\right)^{a}-\frac{b-1}{a-b}\left(Z_{j}\right)^{b}\right]\\&
V_{j}=\left(1-\alpha\right)
\end{split}
\end{equation}

\subsection{Scale invariant F Divergence.}
For such a divergence, the nominal invariance factor has no explicit expression, so we are led to use (for example) the particular invariance factor $K^{*}=\frac{\sum_{j}p_{j}}{\sum_{j}q_{j}}$.\\
This gives the invariant divergence:
\begin{equation}
FI\left(p\|q\right)=\sum_{j}p_{j}\sum_{i}\overline{p}_{i}\log\frac{\overline{p}_{i}}{\alpha \overline{p}_{i}+\left(1-\alpha\right)\overline{q}_{i}}
\end{equation}
With the normalized variables $\overline{p}_{i}=\frac{p_{i}}{\sum_{j}p_{j}}$ and $\overline{q}_{i}=\frac{q_{i}}{\sum_{j}q_{j}}$.\\
Except for the multiplicative factor $\left(\sum_{j}p_{j}\right)$, we recover the expression of the "$F$" initial divergence which is simplified taking into account the normalization of the variables.\\
We will note:
\begin{equation}
\overline{Z}_{i}=\frac{\overline{p}_{i}}{\alpha \overline{p}_{i}+\left(1-\alpha\right)\overline{q}_{i}}
\end{equation} 
By introducing the deformed logarithm, it comes:
\begin{equation}
L_{d}FI\left(p\|q\right)=\sum_{j}p_{j}\sum_{i}\frac{\overline{p}_{i}}{a-b}\left[\left(\overline{Z}_{i}\right)^{a-1}-\left(\overline{Z}_{i}\right)^{b-1}\right]
\end{equation} 
The opposite of the gradient with respect to "$q$" is written $\forall j$:
\begin{equation}
-\frac{\partial L_{d}FI\left(p\|q\right)}{\partial q_{j}}=-\sum_{j}p_{j}\left[\frac{a-1}{a-b}\sum_{i}\overline{p}_{i}\left(\overline{Z}_{i}\right)^{a-2}-\frac{b-1}{a-b}\sum_{i}\overline{p}_{i}\left(\overline{Z}_{i}\right)^{b-2}\right]\frac{\partial \overline{Z}_{i}}{\partial q_{j}} 
\end{equation}
With:
\begin{equation}
\frac{\partial \overline{Z}_{i}}{\partial q_{j}}=\frac{\partial \overline{Z}_{i}}{\partial \overline{q}_{i}}\frac{\partial \overline{q}_{i}}{\partial q_{j}}=-\left[\frac{\left(1-\alpha\right)\overline{p}_{i}}{\left(\alpha \overline{p}_{i}+\left(1-\alpha\right)\overline{q}_{i}\right)^{2}}\right]\left[\frac{\delta_{ij}-\overline{q}_{i}}{\sum_{j}q_{j}}\right] 
\end{equation}
And finally, the expression of the opposite of the gradient with respect to "$q$" is written $\forall j$:
\begin{equation}
\begin{split}
-\frac{\partial L_{d}FI\left(p\|q\right)}{\partial q_{j}}=&-\left(1-\alpha\right)\frac{\sum_{j}p_{j}}{\sum_{j}q_{j}}\left(\frac{a-1}{a-b}\right)\left[\sum_{i}\left(\overline{Z}_{i}\right)^{a}\overline{q}_{i}-\overline{Z}_{j}^{a}\right] 
\\&+\left(1-\alpha\right)\frac{\sum_{j}p_{j}}{\sum_{j}q_{j}}\left(\frac{b-1}{a-b}\right)\left[\sum_{i}\left(\overline{Z}_{i}\right)^{b}\overline{q}_{i}-\overline{Z}_{j}^{b}\right]
\end{split}
\end{equation}
The decomposition of the opposite of the gradient in the form $U-V$, with $U>0$ and $V>0$ is written:
\begin{equation}
\begin{split}
&U_{j}=\left(1-\alpha\right)\frac{\sum_{j}p_{j}}{\sum_{j}q_{j}}\left[\left(\frac{a-1}{a-b}\right)\overline{Z}_{j}^{a}-\left(\frac{b-1}{a-b}\right)\overline{Z}_{j}^{b}\right]>0\\&
V_{j}=\left(1-\alpha\right)\frac{\sum_{j}p_{j}}{\sum_{j}q_{j}}\left[\left(\frac{a-1}{a-b}\right)\sum_{i}\overline{Z}_{i}^{a}\overline{q}_{i}
-\left(\frac{b-1}{a-b}\right)\sum_{i}\overline{Z}_{i}^{b}\overline{q}_{i}\right]>0
\end{split}
\end{equation}

\subsection{G Divergence.}
As previously mentionned in \cite{taneja2001}, \cite{lanteri2019}, \cite{lanteri2020}, it is written:
\begin{equation}
G\left(p\|q\right)=\sum_{i}\left[\alpha p_{i}+\left(1-\alpha\right)q_{i}\right] \log\frac{\alpha p_{i}+\left(1-\alpha\right)q_{i}}{p_{i}}+\left(1-\alpha\right)\left(p_{i}-q_{i}\right) 
\end{equation}
By writting:
\begin{equation}
T_{i}=\alpha p_{i}+\left(1-\alpha\right)q_{i}
\end{equation}
Consequently:
\begin{equation}
\frac{\partial T_{j}}{\partial q_{j}}=1-\alpha
\end{equation}
By introducing the deformed logarithm, it comes:
\begin{equation}
L_{d}G\left(p\|q\right)=\sum_{i}\frac{T_{i}}{a-b}\left[\left(\frac{T_{i}}{p_{i}}\right)^{a-1}-\left(\frac{T_{i}}{p_{i}}\right)^{b-1}\right]+\left(1-\alpha\right)\left(p_{i}-q_{i}\right)  
\end{equation}
After calculations, the expression of the opposite of the gradient with respect to "$q$" is written $\forall j$:
\begin{equation}
-\frac{\partial L_{d}G\left(p\|q\right)}{\partial q_{j}}=\left(1-\alpha\right)\left\lbrace 1-\left[\left(\frac{a}{a-b}\right)\left(\frac{T_{j}}{p_{j}}\right)^{a-1}-\left(\frac{b}{a-b}\right)\left(\frac{T_{j}}{p_{j}}\right)^{b-1}\right]\right\rbrace  
\end{equation}
The decomposition of the opposite of the gradient in the form $U-V$, with $U>0$ and $V>0$ will depend on the sign of $(a-b)$.\\
Indeed, we will have:\\
\textbf{If $a-b>0$:}
\begin{equation}
U_{j}=(1-\alpha)\left[1+\left(\frac{b}{a-b}\right)\left(\frac{T_{j}}{p_{j}}\right)^{b-1}\right]\;\;;\;\;V_{j}=\left(\frac{a}{a-b}\right)\left(\frac{T_{j}}{p_{j}}\right)^{a-1}
\end{equation}
\textbf{If $a-b<0$:}
\begin{equation}
U_{j}=(1-\alpha)\left[1+\left(\frac{a}{b-a}\right)\left(\frac{T_{j}}{p_{j}}\right)^{a-1}\right]\;\;;\;\;V_{j}=\left(\frac{b}{b-a}\right)\left(\frac{T_{j}}{p_{j}}\right)^{b-1}
\end{equation}

\subsection{Scale invariant G Divergence.}
For such divergence, the nominal invariance factor has no explicit expression, so we are led to use (for example) the particular invariance factor $K^{*}=\frac{\sum_{j}p_{j}}{\sum_{j}q_{j}}$.\\
This gives the invariant divergence:
\begin{equation}
GI\left(p\|q\right)=\sum_{j}p_{j}\sum_{i}\left[\alpha \overline{p}_{i}+\left(1-\alpha\right)\overline{q}_{i}\right] \log\frac{\alpha \overline{p}_{i}+\left(1-\alpha\right)\overline{q}_{i}}{\overline{p}_{i}}
\end{equation}
With the normalized variables $\overline{p}_{i}=\frac{p_{i}}{\sum_{j}p_{j}}$ and $\overline{q}_{i}=\frac{q_{i}}{\sum_{j}q_{j}}$.\\
Except for the multiplicative factor $\left(\sum_{j}p_{j}\right)$, we recover the expression of the "$G$" initial divergence.
By writting:
\begin{equation}
\overline{T}_{i}=\alpha \overline{p}_{i}+\left(1-\alpha\right)\overline{q}_{i}
\end{equation}
By introducing the deformed logarithm, it comes:
\begin{equation}
L_{d}GI\left(p\|q\right)=\frac{\sum_{j}p_{j}}{a-b}\sum_{i}\left[\frac{\overline{T}^{a}_{i}}{\overline{p}^{a-1}_{i}}-\frac{\overline{T}^{b}_{i}}{\overline{p}^{b-1}_{i}}\right] 
\end{equation}
The opposite of the gradient with respect to "$q$" is written $\forall j$:
\begin{equation}
-\frac{\partial L_{d}GI\left(p\|q\right)}{\partial q_{j}}=-\sum_{j}p_{j}\left[\frac{a}{a-b}\sum_{i}\left(\frac{\overline{T}_{i}}{\overline{p}_{i}}\right)^{a-1}-\frac{b}{a-b}\sum_{i}\left(\frac{\overline{T}_{i}}{\overline{p}_{i}}\right)^{b-1}\right]\frac{\partial \overline{T}_{i}}{\partial q_{j}} 
\end{equation}
With:
\begin{equation}
\frac{\partial \overline{T}_{i}}{\partial q_{j}}=\frac{\partial \overline{T}_{i}}{\partial \overline{q}_{i}}\frac{\partial \overline{q}_{i}}{\partial q_{j}}=\left(1-\alpha\right)\frac{\delta_{ij}-\overline{q}_{i}}{\sum_{j}q_{j}}
\end{equation}
Hence, finally:
\begin{equation}
\begin{split}
-\frac{\partial L_{d}GI\left(p\|q\right)}{\partial q_{j}}=&-\left(1-\alpha\right)\frac{\sum_{j}p_{j}}{\sum_{j}q_{j}}\left[\frac{a}{a-b}\left(\frac{\overline{T}_{j}}{\overline{p}_{j}}\right)^{a-1}-\frac{a}{a-b}\sum_{i}\left(\frac{\overline{T}_{i}}{\overline{p}_{i}}\right)^{a-1}\overline{q}_{i}\right]\\&+\left(1-\alpha\right)\frac{\sum_{j}p_{j}}{\sum_{j}q_{j}}\left[\frac{b}{a-b}\left(\frac{\overline{T}_{j}}{\overline{p}_{j}}\right)^{b-1}-\frac{b}{a-b}\sum_{i}\left(\frac{\overline{T}_{i}}{\overline{p}_{i}}\right)^{b-1}\overline{q}_{i}\right] 
\end{split}
\end{equation}
The decomposition of the opposite of the gradient in the form $U-V$, with $U>0$ and $V>0$ will depend on the sign of $(a-b)$.\\
Indeed, we will have:\\
\textbf{If $a-b>0$:}
\begin{equation}
\begin{split}
&U_{j}=(1-\alpha)\frac{\sum_{j}p_{j}}{\sum_{j}q_{j}}\left[\frac{a}{a-b}\sum_{i}\left(\frac{\overline{T}_{i}}{\overline{p}_{i}}\right)^{a-1}\overline{q}_{i}+\frac{b}{a-b}\left(\frac{\overline{T}_{j}}{\overline{p}_{j}}\right)^{b-1}\right]\\&V_{j}=(1-\alpha)\frac{\sum_{j}p_{j}}{\sum_{j}q_{j}}\left[\frac{b}{a-b}\sum_{i}\left(\frac{\overline{T}_{i}}{\overline{p}_{i}}\right)^{b-1}\overline{q}_{i}+\frac{a}{a-b}\left(\frac{\overline{T}_{j}}{\overline{p}_{j}}\right)^{a-1}\right]
\end{split}
\end{equation}

\textbf{If $a-b<0$:}
\begin{equation}
\begin{split}
&U_{j}=(1-\alpha)\frac{\sum_{j}p_{j}}{\sum_{j}q_{j}}\left[\frac{b}{b-a}\sum_{i}\left(\frac{\overline{T}_{i}}{\overline{p}_{i}}\right)^{b-1}\overline{q}_{i}+\frac{a}{b-a}\left(\frac{\overline{T}_{j}}{\overline{p}_{j}}\right)^{a-1}\right]\\&V_{j}=(1-\alpha)\frac{\sum_{j}p_{j}}{\sum_{j}q_{j}}\left[\frac{a}{b-a}\sum_{i}\left(\frac{\overline{T}_{i}}{\overline{p}_{i}}\right)^{a-1}\overline{q}_{i}+\frac{b}{b-a}\left(\frac{\overline{T}_{j}}{\overline{p}_{j}}\right)^{b-1}\right]
\end{split}
\end{equation}

\section{Algorithmics.}
\subsection{Overview of the construction method of the algorithms.}
We recall that in the case that we are dealing with (linear model), we have $p=y$, $q=Hx$, thus $q_{j}^{k}=(Hx^{k})_{j}$.\\
Consequently, for a divergence $D(p\|q)$ or $D(q\|p)$, we have:
\begin{equation}
\frac{\partial D}{\partial x}=H^{T}\frac{\partial D}{\partial q}	
\end{equation}
This justifies the fact that throughout this work, we have developed the expressions of the gradients with respect to ``$q$''.\\
The algorithms proposed here are based on the SGM method or its variants, as described in detail in \cite{lanteri2019} and \cite{lanteri2020}.\
If the considered divergence is of classical non-invariant form, the proposed algorithms take into account the non-negativity constraint of the solution.\\
On the other hand, if we want to take into account a sum constraint on the components of the solution, the scale invariant divergences must be used.\\
In any case, the general form of the algorithms remains the same.\\
The basic iterative algorithm is written in general form:
\begin{equation}
x^{k+1}_{j}=x^{k}_{j}+\alpha^{k}_{j}x^{k}_{j}\left(-\frac{\partial D}{\partial x}\right)^{k}_{j}
\label{algogen}	
\end{equation}
This writing is always possible.\\              
The divergence $D$ being convex with respect to ``$x$'', $\left(-\frac{\partial D}{\partial x}\right)$ is a direction of descent and the vector of components $\left[x_{j}.\left(-\frac{\partial D}{\partial x}\right)_{j}\right] $ is a direction of descent if $x_{j}>0$.\\
Consequently, with a positive initial estimate, as long as the operating mode ensures the non-negativity of the successive iterates, such an algorithm converges provided that the descent step is properly computed; these two points will be specified in what follows.\\
To our knowledge, with the only notable exception of the dual Kullback-Leibler divergence which will be developed in Appendix 3, the opposite of the gradient can always be written in the form:
\begin{equation}
\left(-\frac{\partial D}{\partial x}\right)^{k}_{j}=U^{k}_{j}-V^{k}_{j}\ \ ;\ \ \ \ U^{k}_{j}>0\ \ ;\ \ \ \ V^{k}_{j}>0	
\end{equation}
In some (very rare) cases, this decomposition requires some reflection.\\ 
Hence the almost always possible writing:
\begin{equation}
x^{k+1}_{j}=x^{k}_{j}+\alpha^{k}_{j}x^{k}_{j}\left(U^{k}_{j}-V^{k}_{j}\right)	
\end{equation}
Then, preconditioning when possible by $\frac{1}{V^{k}_{j}}>0$, we obtain the pseudo-multiplicative form:
\begin{equation}
x^{k+1}_{j}=x^{k}_{j}+\alpha^{k}_{j}x^{k}_{j}\left(\frac{U^{k}_{j}}{V^{k}_{j}}-1\right)	
\end{equation}
Due to $V^{k}_{j}>0$ $\forall j$, the opposite of the modified gradient remains a descent direction.\\\\
Regarding the descent stepsize, for each of these two algorithms, the following procedure will be implemented:\\
* 1- At a given iteration ``$k$'', we compute the maximum step $(\alpha^{k})_{Max}$ ensuring the non-negativity of all components of $x^{k+1}$.\\
* 2 - The descent stepsize $\alpha^{k}$ (valid for all components), ensuring the convergence of the algorithm is then computed by a one-dimensional search method such as Armijo \cite{armijo1966} (for example), over the interval $\left[0,(\alpha^{k})_{Max}\right]$.\\
This procedure is analyzed in a more detailed way in \cite{lanteri2019} and \cite{lanteri2020}.\\
We can write in a general way, whatever the expression of the opposite of the gradient, with a descent step independent of the component:
\begin{equation}
x^{k+1}=x^{k}+\alpha^{k}x^{k}\left(U^{k}-V^{k}\right)
\label{algobase}	
\end{equation}

\textbf{Remark 1}: in this expression, the operation $x^{k}\left(U^{k}-V^{k}\right)$ represents the component to component product of the vectors $x^{k}$ and $\left(U^{k}-V^{k}\right)$ (Hadamard product).\\
With a modified (preconditioned) gradient, we will obtain the pseudo-multiplicative form:
\begin{equation}
x^{k+1}=x^{k}+\alpha^{k}x^{k}\left(\frac{U^{k}}{V^{k}}-1\right)
\label{algoprecond}	
\end{equation}
For this expression, the previous remark applies.\\

\textbf{Remark 2}: in this notation, $\frac{U^{k}}{V^{k}}$ is a vector obtained by making the ratio component by component of the vectors $U^{^{k}}$ and $V^{^{k}}$, as well as ``$1$'' is the unit vector.\\ 
In the latter case, if we use a descent step $\alpha^{k}=1\ ,\forall k$, the non-negativity is ensured, and we obtain a purely multiplicative algorithm which is written:
\begin{equation}
x^{k+1}=x^{k}\left(\frac{U^{k}}{V^{k}}\right)
\label{algomult}	
\end{equation}
Of course, in all generality, nothing proves the convergence of purely multiplicative algorithms, each divergence implies a particular analysis. \\

\textbf{Remark 3} : For non-invariant divergences, the algorithms (\ref{algobase}) (\ref{algoprecond}) and (\ref{algomult}) only ensure the non-negativity of the solution.\\

\textbf{Remark 4} : The scale invariant divergences make all their sense if one requires in addition that the sum constraint is ensured, indeed, an algorithm of type (\ref{algobase}) applied on such divergences allows to ensure the property:
\begin{equation}
	\sum_{l}x^{k+1}_{l}=\sum_{l}x^{k}_{l}
\end{equation}
Starting from an initial estimate $x^{0}$ such that $sum_{l}x^{0}_{l}=C$, all successive estimates will be of the same sum.\\
However, when such scale invariant divergences are considered, the use of a pseudo multiplicative (preconditioned) algorithm (\ref{algoprecond}) or of a purely multiplicative algorithm of the type (\ref{algomult}) (as long as its convergence is ensured), does not allow to ensure spontaneously the sum constraint; an additional step is necessary.\\
So, at each iteration, the procedure is carried out in 2 steps:\\\\
* - We first calculate a preliminary estimate:
\begin{equation}
\widetilde{x}^{k+1}=x^{k}+\alpha^{k}x^{k}\left(\frac{U^{k}}{V^{k}}-1\right)
\label{estprovbis}	
\end{equation}
or:
\begin{equation}
\widetilde{x}^{k+1}=x^{k}\left(\frac{U^{k}}{V^{k}}\right)
\label{estprov}	
\end{equation}

* - Then in a normalization step, we compute:
\begin{equation}
	x^{k+1}=\frac{\widetilde{x}^{k+1}}{\sum_{l}\tilde{x}^{k+1}_{l}}C
\label{estfinal}	
\end{equation}
Given the properties of the scale invariant divergences, this last operation does not change the value of the divergence under consideration.

\section{Appendix 1}
The deformed logarithm applied to the product $(x.y)$ is written:
\begin{equation}
	L_{d}\left(x.y\right)=\frac{(x.y)^{a-1}-(x.y)^{b-1}}{a-b}
\end{equation}
With for example $a=1+\epsilon$ and $b=1-\epsilon$, we obtain the simplified form:
\begin{equation}
LS_{d}\left(x.y\right)=\frac{(x.y)^{\epsilon}-(x.y)^{-\epsilon}}{2\;\epsilon}	
\end{equation}
Which can be written:
\begin{equation}
LS_{d}\left(x.y\right)=\frac{\exp \epsilon \log (x.y)-\exp (-\epsilon) \log (x.y)}{2\;\epsilon}	
\end{equation}
Taking the first order limited expansion, it comes:
\begin{equation}
LS_{d}\left(x.y\right)=\frac{1+\epsilon \log (x.y)-(1-\epsilon \log (x.y))}{2\;\epsilon}	
\end{equation}
Then:
\begin{equation}
LS_{d}\left(x.y\right)=\log (x.y)=\log x+\log y	
\end{equation}
Q.E.D.

\section{Appendix 2}
The generalization of divergences between means proposed by Taneja \cite{taneja1989} can be summarized as follows.\\
For instance, let us consider the inequality between the weighted arithmetic mean $\left(MA\right)_{i}=\alpha p_{i}+\left(1-\alpha\right) q_{i}$ and the weighted geometric mean $\left(MG\right)_{i}=p^{\alpha}_{i}q^{1-\alpha}_{i}$; 
\begin{equation}
\left(MG\right)_{i}-\left(MA\right)_{i}\leq0\ \ \ \Rightarrow\ \ \left(MG\right)_{i}^{1-r}-\left(MA\right)_{i}^{1-r}\leq0	\ \ ; \ \ r\leq1
\end{equation}
And,then:
\begin{equation}
\left(MG\right)_{i}^{1-r}\left(MA\right)_{i}^{r}-\left(MA\right)_{i}\leq0	
\end{equation}
The extension to the whole field leads to:
\begin{equation}
\sum_{i}\left(MG\right)_{i}^{1-r}\left(MA\right)_{i}^{r}-\sum_{i}\left(MA\right)_{i}\leq0	
\end{equation}
And finally:
\begin{equation}
\frac{1}{r-1}\left[\sum_{i}\left(MG\right)_{i}^{1-r}\left(MA\right)_{i}^{r}-\sum_{i}\left(MA\right)_{i}\right]\geq0	
\end{equation}
From this point, we can apply on each of the terms of the difference, an increasing function without changing the sign of the inequality.\\
Here, Tanéja uses the Tsallis deformed logarithm (\ref{logdtsa}) with the parameter $\left(s-1/r-1\right)$ instead of $t-1$, to obtain the divergence:
\begin{equation}
D^{s}_{r}\left(p\|q\right)=\frac{1}{s-1}
\left\{\left[\sum_{i}\left(MG\right)_{i}^{1-r}\left(MA\right)_{i}^{r}\right]^{\frac{s-1}{r-1}}-\left[\sum_{i}\left(MA\right)_{i}\right]^{\frac{s-1}{r-1}}\right\}	
\end{equation}
In fact, Taneja proposes something less general than this, because he limits himself to the case of variables of sum equal to "1", in this case, the 2nd term of the difference will be equal to "1".\\

\section{Appendix 3}
\subsection{The special case of the dual Kullback-Leibler divergence.}

It is written:
\begin{equation}
KL(q\|p)=\sum_{i}q_{i}\log\frac{q_{i}}{p_{i}}+p_{i}-q_{i}
\label{divklduale}
\end{equation}
The opposite of its gradient with respect to "$q$" is written $\forall j$:
\begin{equation}
-\frac{\partial KL(q\|p)}{\partial q_{j}}=\log p_{j}-\log q_{j}
\end{equation}
This is where the problem arises, indeed, a decomposition of this expression in the form $U-V$ with the strict inequalities $U>0$ and $V>0$ is impossible.\\
This decomposition would be:\\

* If $q_{j}>1$ et $p_{j}<1$ , the opposite of its gradient is negative,\\
 $U_{j}=0\;\;;\;\;V_{j}=-\log p_{j}+\log q_{j}$.\\
 
 * If $q_{j}<1$ et $p_{j}>1$ , the opposite of its gradient is positive,\\
 $U_{j}=\log p_{j}-\log q_{j}\;\;;\;\;V_{j}=0$.\\
 
 * If $q_{j}>1$ et $p_{j}>1$,\\
 $U_{j}=\log p_{j}\;\;;\;\;V_{j}=\log q_{j}$.\\
 
 * If $q_{j}<1$ et $p_{j}<1$,\\
  $U_{j}=-\log q_{j}\;\;;\;\;V_{j}=-\log p_{j}$.\\
  
 Thus, if an algorithm of the general form (\ref{algogen}) can always be used, it is not the same for algorithms involving a decomposition of the opposite of the gradient, such as the pseudo- multiplicative (\ref{algoprecond}) or purely multiplicative (\ref{algomult}) algorithms.\\ 
Such a problem does not arise (apparently) when we consider the form of this divergence obtained by introducing the deformed logarithm.\\
Indeed by replacing the natural logarithm by the deformed logarithm in the expression (\ref{divklduale}), we obtain:
\begin{equation}
L_{d}KL(q\|p)=\frac{1}{a-b}\sum_{i}\left[\frac{q_{i}^{a}}{p_{i}^{a-1}}-\frac{q_{i}^{b}}{p_{i}^{b-1}}\right]+\sum_{i}p_{i}-\sum_{i}q_{i}
\end{equation} 
The opposite of its gradient with respect to "$q$" is written $\forall j$:
\begin{equation}
-\frac{\partial L_{d}KL(q\|p)}{\partial q_{j}}=1-\left[\frac{a}{\left(a-b\right)}\frac{q_{i}^{a-1}}{p_{i}^{a-1}}-\frac{b}{\left(a-b\right)}\frac{q_{i}^{b-1}}{p_{i}^{b-1}}\right] 
\end{equation}    
The decomposition of the opposite of the gradient in the form $U-V$ with $U>0$ and $V>0$ depends on the sign of $(a-b)$, it is written:\\

\textbf{* If $a-b>0$}
\begin{equation}
U_{j}=1+\frac{b}{a-b}\left(\frac{q_{j}}{p_{j}}\right)^{b-1}\;\;\;;\;\;\;V_{j}=\frac{a}{a-b}\left(\frac{q_{j}}{p_{j}}\right)^{a-1} 
\end{equation}

\textbf{* If $a-b<0$}
\begin{equation}
U_{j}=1+\frac{a}{b-a}\left(\frac{q_{j}}{p_{j}}\right)^{a-1}\;\;\;;\;\;\;V_{j}=\frac{b}{b-a}\left(\frac{q_{j}}{p_{j}}\right)^{b-1} 
\end{equation}  
  
 \subsection{Scale invariant form.}
 The problem above mentioned also appears for this divergence.\\
The nominal invariance factor is explicitly computable \cite{lanteri2019}, \cite{lanteri2020}.\\
Its expression is:
\begin{equation}
K_{0}\propto \exp \left(\frac{\sum_{i}q_{i}\log \frac{p_{i}}{q_{i}} }{\sum_{j}q_{j}}\right) 
\end{equation}
The corresponding invariant divergence is written:
\begin{equation}
KLI_{n}(q\|p)=\sum_{i}p_{i}-K_{0}\sum_{i}q_{i}
\end{equation}
The opposite of the gradient with respect to "$q$" is written, $forall j$:
\begin{equation}
-\frac{\partial KLI_{n}(q\|p)}{\partial q_{j}}=K_{0}\left[\log \frac{p_{j}}{q_{j}}-\frac{\sum_{i}q_{i}\log \frac{p_{i}}{q_{i}} }{\sum_{j}q_{j}}\right] 
\end{equation}
Thus, we can observe that with such an invariance factor, a decomposition of this expression in the form $U-V$ with $U>0$ and $V>0$ encounters the same problems as those already pointed out in the previous sub-section for the non-invariant divergence, with the same algorithmic problems.\\
For comparison, we develop the case in which we use the invariance factor 
 $K^{*}=\frac{\sum_{j}p_{j}}{\sum_{j}q_{j}}$.\\
 With $\overline{p}_{i}=\frac{p_{i}}{\sum_{j}p_{j}}$ and $\overline{q}_{i}=\frac{q_{i}}{\sum_{j}q_{j}}$, the invariant divergence thus obtained is written after simplification: 
\begin{equation}
KLI(q\|p)=\sum_{j}p_{j}\sum_{i}\overline{q}_{i}\log \frac{\overline{q}_{i}}{\overline{p}_{i}}
\label{kldualeinv}
\end{equation}
The opposite of the gradient with respect to "$q$" is written all calculations done, $\forall j$:
\begin{equation}
-\frac{\partial KL(q\|p)}{\partial q_{j}}=\frac{\sum_{j}p_{j}}{\sum_{j}q_{j}}\left(\sum_{i}\overline{q}_{i}\log \frac{\overline{q}_{i}}{\overline{p}_{i}}-\log \frac{\overline{q}_{j}}{\overline{p}_{j}}\right) 
\end{equation}
In this case, the decomposition of the opposite of the gradient as mentioned before is no more a problem due to the normalization of the variables.\\
Indeed, we will have $\forall j$:
\begin{equation}
U_{j}=-\sum_{i}\overline{q}_{i}\log \overline{p}_{i}-\log \overline{q}_{j}\;\;;\;\;V_{j}=-\sum_{i}\overline{q}_{i}\log \overline{q}_{i}-\log \overline{p}_{j}
\end{equation}

We show, as an example, that the problem of gradient decomposition do not appear when we introduce the deformed logarithm in the expression (\ref{kldualeinv}).\\
We then have:
\begin{equation}
L_{d}KLI(q\|p)=\sum_{j}p_{j}\sum_{i}\frac{\overline{q}_{i}}{a-b}\left[\left(\frac{\overline{q}_{i}}{\overline{p}_{i}}\right)^{a-1}-\left(\frac{\overline{q}_{i}}{\overline{p}_{i}}\right)^{b-1}\right]  
\end{equation}
The opposite of the gradient with respect to "$q$" is written, $\forall j$:
\begin{equation}
-\frac{\partial L_{d}KLI(q\|p)}{\partial q_{j}}=-\sum_{j}p_{j}\sum_{i}\left[\frac{a}{a-b}\left( \frac{\overline{q}_{i}}{\overline{p}_{i}}\right)^{a-1}-\frac{b}{a-b}\left( \frac{\overline{q}_{i}}{\overline{p}_{i}}\right)^{b-1}\right]\frac{\partial \overline{q}_{i}}{\partial q_{j}}
\end{equation}
This leads, after all calculations, to:
\begin{equation}
\begin{split}
-\frac{\partial L_{d}KLI(q\|p)}{\partial q_{j}}=&\frac{\sum_{j}p_{j}}{\sum_{j}q_{j}}\left[ \sum_{i}\frac{a}{a-b}\left( \frac{\overline{q}_{i}}{\overline{p}_{i}}\right)^{a}-\sum_{i}\frac{b}{a-b}\left( \frac{\overline{q}_{i}}{\overline{p}_{i}}\right)^{b}\right] 
\\&-\frac{\sum_{j}p_{j}}{\sum_{j}q_{j}}\left[\frac{a}{a-b}\left( \frac{\overline{q}_{j}}{\overline{p}_{j}}\right)^{a-1}-\frac{b}{a-b}\left(\frac{\overline{q}_{j}}{\overline{p}_{j}}\right)^{b-1}\right]  
\end{split}
\end{equation}
With this expression the decomposition in the form $U-V$ with $U>0$ and $V>0$ is not a problem but depends, of course, on the sign of $(a-b)$, and we will have:\\

\textbf{* If $a-b>0$}
\begin{equation}
\begin{split}
&U_{j}=\frac{\sum_{j}p_{j}}{\sum_{j}q_{j}}\left[ \sum_{i}\frac{a}{a-b}\;\;\frac{\overline{q}_{i}^{a}}{\overline{p}_{i}^{a-1}}+\frac{b}{a-b}\left(\frac{\overline{q}_{j}}{\overline{p}_{j}}\right)^{b-1}\right]\\&
V_{j}=\frac{\sum_{j}p_{j}}{\sum_{j}q_{j}}\left[ \sum_{i}\frac{b}{a-b}\;\;\frac{\overline{q}_{i}^{b}}{\overline{p}_{i}^{b-1}}+\frac{a}{a-b}\left(\frac{\overline{q}_{j}}{\overline{p}_{j}}\right)^{a-1}\right]
\end{split}
\end{equation}

\textbf{* If $a-b<0$}:
\begin{equation}
\begin{split}
&U_{j}=\frac{\sum_{j}p_{j}}{\sum_{j}q_{j}}\left[ \sum_{i}\frac{b}{b-a}\;\;\frac{\overline{q}_{i}^{b}}{\overline{p}_{i}^{b-1}}+\frac{a}{b-a}\left(\frac{\overline{q}_{j}}{\overline{p}_{j}}\right)^{a-1}\right]\\&
V_{j}=\frac{\sum_{j}p_{j}}{\sum_{j}q_{j}}\left[ \sum_{i}\frac{a}{b-a}\;\;\frac{\overline{q}_{i}^{a}}{\overline{p}_{i}^{a-1}}+\frac{b}{b-a}\left(\frac{\overline{q}_{j}}{\overline{p}_{j}}\right)^{b-1}\right]
\end{split}
\end{equation}

\bibliographystyle{plain}
\bibliographystyle{alpha}
\bibliography{biblio}

\begin{thebibliography}{10}

\bibitem{abe1997}
S.~Abe.
\newblock A note on the q-deformation-theoretic aspect of the generalized
  entropies in nonextensive physics.
\newblock {\em Physics Letters A}, 224(6):326--330, 1997.

\bibitem{amari2009}
S.I. Amari.
\newblock Alpha divergence is unique, belonging to both $f$-divergence and
  {B}regman divergence classes.
\newblock {\em Information Theory, IEEE Transactions on}, 55(11):4925--4931,
  2009.

\bibitem{armijo1966}
L.~Armijo.
\newblock Minimization of functions having {L}ipschitz continuous first partial
  derivatives.
\newblock {\em Pacific Journal of mathematics}, 16(1):1--3, 1966.

\bibitem{ben1989}
A.~Ben-Tal, A.~Charnes, and M.~Teboulle.
\newblock Entropic means.
\newblock {\em Journal of Mathematical Analysis and Applications},
  139(2):537--551, 1989.

\bibitem{bertero1998}
M.~Bertero and P.~Boccacci.
\newblock {\em Introduction to inverse problems in imaging}.
\newblock CRC press, 1998.

\bibitem{borges1998}
E.P. Borges and I.~Roditi.
\newblock A family of nonextensive entropies.
\newblock {\em Physics Letters A}, 246(5):399--402, 1998.

\bibitem{chakrabarti1991}
R.~Chakrabarti and R.~Jagannathan.
\newblock A (p, q)-oscillator realization of two-parameter quantum algebras.
\newblock {\em Journal of Physics A: Mathematical and General}, 24(13):L711,
  1991.

\bibitem{cichocki2010}
A.~Cichocki and S.I. Amari.
\newblock Families of $\alpha$-$\beta$-and $\gamma$-divergences: Flexible and
  robust measures of similarities.
\newblock {\em Entropy}, 12(6):1532--1568, 2010.

\bibitem{cichocki2011}
A.~Cichocki, S.~Cruces, and S.I. Amari.
\newblock Generalized alpha-beta divergences and their application to robust
  nonnegative matrix factorization.
\newblock {\em Entropy}, 13(1):134--170, 2011.

\bibitem{fevotte2011}
C.~F{\'e}votte and J.~Idier.
\newblock Algorithms for nonnegative matrix factorization with the
  $\beta$-divergence.
\newblock {\em Neural computation}, 23(9):2421--2456, 2011.

\bibitem{furuichi2010}
S.~Furuichi.
\newblock An axiomatic characterization of a two-parameter extended relative
  entropy.
\newblock {\em Journal of mathematical physics}, 51(12):123302, 2010.

\bibitem{Idier01a}
J.~Idier, editor.
\newblock {\em Approche bay\'esienne pour les probl\`emes inverses}.
\newblock Trait\'e IC2, S\'erie traitement du signal et de l'image, Herm\`es,
  Paris, nov. 2001.

\bibitem{jackson1909}
F.H. Jackson.
\newblock q-form of taylor's theorem.
\newblock {\em Mess. Math}, 3:57, 1909.

\bibitem{jackson1910}
F.H. Jackson.
\newblock On q—definite integrals, j.
\newblock {\em Pure Appl. Math}, 41:193--403, 1910.

\bibitem{kaniadakis2002}
G.~Kaniadakis.
\newblock Statistical mechanics in the context of special relativity.
\newblock {\em Physical review E}, 66(5):056125, 2002.

\bibitem{kaniadakis2005}
G.~Kaniadakis, M.~Lissia, and A.M. Scarfone.
\newblock Two-parameter deformations of logarithm, exponential, and entropy: A
  consistent framework for generalized statistical mechanics.
\newblock {\em Physical Review E}, 71(4):046128, 2005.

\bibitem{kullback1951}
S.~Kullback and R.A. Leibler.
\newblock On information and sufficiency.
\newblock {\em The annals of mathematical statistics}, 22(1):79--86, 1951.

\bibitem{lanteri2019}
H.~Lant{\'e}ri.
\newblock Divergences. divergences invariantes. applications aux probl{\`e}mes
  inverses lin{\'e}aires, nmf et d{\'e}convolution aveugle.
\newblock {\em HAL.Archives ouvertes hal-01745256v3}, 2019.

\bibitem{lanteri2020}
H.~Lant{\'e}ri.
\newblock Divergences. scale invariant divergences. applications to linear
  inverse problems. nmf blind deconvolution.
\newblock {\em arXiv preprint arXiv:2003.01411}, 2020.

\bibitem{lanteri2021}
H.~Lant{\'e}ri.
\newblock Deformed logarithms. associated entropic divergences. applications to
  linear inverse problems. inversion algorithms.
\newblock {\em arXiv preprint arXiv:2109.12895}, 2021.

\bibitem{lanteri2022}
H.~Lant{\'e}ri.
\newblock Logarithmes déformés. -divergences entropiques associées.
  -applications aux problèmes inverses linéaires. - algorithmes d'inversion.
\newblock {\em HAL.Archives ouvertes hal-03355082v3}, 2022.

\bibitem{mcanally1995}
D.S. McAnally.
\newblock q-exponential and q-gamma functions. i. q-exponential functionsa.
\newblock {\em Journal of Mathematical Physics}, 36(1):546--573, 1995.

\bibitem{mittal1975}
D.P. Mittal.
\newblock On some functional equations concerning entropy, directed divergence
  and inaccuracy.
\newblock {\em Metrika}, 22(1):35--45, 1975.

\bibitem{naudts2002}
J.~Naudts.
\newblock Deformed exponentials and logarithms in generalized thermostatistics.
\newblock {\em Physica A: Statistical Mechanics and its Applications},
  316(1-4):323--334, 2002.

\bibitem{shannon1948}
C.~E. Shannon.
\newblock A mathematical theory of communication.
\newblock {\em Bell system technical journal}, 27(3):379--423, 1948.

\bibitem{sharma1975}
B.D. Sharma and I.J. Taneja.
\newblock Entropy of type ($\alpha$, $\beta$) and other generalized measures in
  information theory.
\newblock {\em Metrika}, 22(1):205--215, 1975.

\bibitem{taneja2001}
I.J. Taneja.
\newblock Generalized information measures and their applications. on-line
  book, 2001.
\newblock {\em URL www. mtm. ufsc. br/taneja/book/book. html}.

\bibitem{taneja1989}
I.J. Taneja.
\newblock On generalized information measures and their applications.
\newblock {\em Advances in Electronics and Electron Physics}, 76:327--413,
  1989.

\bibitem{tsallis1988}
C.~Tsallis.
\newblock Possible generalization of boltzmann-gibbs statistics.
\newblock {\em Journal of statistical physics}, 52(1):479--487, 1988.

\bibitem{wada2010}
T.~Wada and A.M. Scarfone.
\newblock Finite difference and averaging operators in generalized entropies.
\newblock In {\em Journal of Physics: Conference Series}, volume 201, page
  012005. IOP Publishing, 2010.

\end{thebibliography}

\end{document}